\documentclass[review,onefignum,onetabnum]{siamart171218}

\usepackage{lipsum}
\usepackage{amsfonts}
\usepackage{graphicx}
\usepackage{epstopdf}
\usepackage{algorithmic}
\usepackage{amsmath,amssymb}
\usepackage{textcomp,marvosym}
\usepackage{cite}
\usepackage{nameref,hyperref}
\usepackage{changes}
\usepackage{subcaption}
\usepackage{float}
\usepackage[T1]{fontenc}
\ifpdf
  \DeclareGraphicsExtensions{.eps,.pdf,.png,.jpg}
\else
  \DeclareGraphicsExtensions{.eps}
\fi

\newcommand{\bs}[1]{\boldsymbol{#1}}
\newsiamremark{remark}{Remark}
\newsiamremark{hypothesis}{Hypothesis}
\crefname{hypothesis}{Hypothesis}{Hypotheses}
\newsiamthm{claim}{Claim}

\headers{Bayesian Calibration of Cell Migration Models}{C. Schenk, J. Ayensa-Jim\'enez, I. Romero}

\title{Bayesian Calibration and Model Assessment of Cell Migration Dynamics with Surrogate Model Integration\thanks{The first two authors are co-first authors. Submitted to the editors 09/23/2025.
\funding{The authors acknowledge financial support from the Spanish Ministry of Science and Innovation through a Juan de la Cierva Postdoctoral Fellowship for J.A-J. (Grant No. JDC2023-052319-I).}}}

\author{Christina Schenk\thanks{IMDEA Materials Institute, C/ Eric Kandel 2, Getafe, 28906, Spain
  (\email{christina.schenk@imdea.org}).}
\and Jacobo Ayensa Jim\'enez\thanks{Universidad Polit\'ecnica de Madrid, Jos\'e Guti\'errez Abascal, 2, Madrid, 29006, Spain,\\TME Lab, Aragón Institute of Engineering Research (I3A), Mariano Esquillor, S/N, Zaragoza, 50018, Spain
  (\email{jacobo.ayensa@upm.es}).}
\and Ignacio Romero\footnotemark[3]\thanks{Universidad Polit\'ecnica de Madrid, Jos\'e Guti\'errez Abascal, 2, Madrid, 29006, Spain\\IMDEA Materials Institute, C/ Eric Kandel 2, Getafe, 28906, Spain (\email{ignacio.romero@upm.es})}}
\usepackage{amsopn}

\makeatletter
\newcommand*{\addFileDependency}[1]{
  \typeout{(#1)}
  \@addtofilelist{#1}
  \IfFileExists{#1}{}{\typeout{No file #1.}}
}
\makeatother

\ifpdf
\hypersetup{
  pdftitle={Bayesian Calibration of Cell Migration Models},
  pdfauthor={C. Schenk, J. Ayensa Jim\'enez, I. Romero}
}
\fi

\nolinenumbers
\begin{document}

\maketitle

\begin{abstract}
Computational models provide crucial insights into complex biological processes such as cancer evolution, but their mechanistic nature often makes them nonlinear and parameter-rich, complicating calibration. We systematically evaluate parameter probability distributions in cell migration models using Bayesian calibration across four complementary strategies: parametric and surrogate models, each with and without explicit model discrepancy. This approach enables joint analysis of parameter uncertainty, predictive performance, and interpretability. Applied to a real data experiment of glioblastoma progression in microfluidic devices, surrogate models achieve higher computational efficiency and predictive accuracy, whereas parametric models yield more reliable parameter estimates due to their mechanistic grounding. Incorporating model discrepancy exposes structural limitations, clarifying where model refinement is necessary. Together, these comparisons offer practical guidance for calibrating and improving computational models of complex biological systems.

\textbf{Relevance to Life Sciences}
Glioblastoma is an aggressive brain tumor where invasive migration drives poor treatment outcomes. Comparing calibration approaches for models of glioblastoma progression in bio-mimetic platforms highlights key trade-offs: surrogate models improve prediction and speed, while parametric models preserve mechanistic interpretability relevant for biological insight. Explicit discrepancy modeling identifies when predictions fail to capture observed tumor behavior, guiding biologically motivated model improvements. These findings demonstrate how computational strategies can advance cancer modeling and support therapeutic planning.

\textbf{Mathematical Content}
We employ Bayesian calibration to estimate posterior parameter distributions across four frameworks: parametric versus surrogate models, each assessed with and without discrepancy terms. Surrogates accelerate inference through approximated dynamics, while parametric formulations directly encode migration processes. The model discrepancy is formalized as an additive term that captures structural gaps between predictions and data. This systematic comparison reveals trade-offs between interpretability, uncertainty quantification, and efficiency, establishing a principled methodology for nonlinear model calibration in life sciences.
\end{abstract}



\section{Introduction}
\label{sec:introduction}
Cell migration models are critical for predicting cancer progression. The complexity, heterogeneity, and dynamic nature of the Tumor Microenvironment hinder \emph{in vivo} investigation of all tumor–stroma interactions. To overcome this, laboratory experiments at multiple scales aim to replicate the spatial and structural features of physiological and pathological environments. However, the limited predictive power of current \emph{in vitro} models remains a major factor in the ongoing decline of new drug approvals relative to research investment \cite{scannell2012diagnosing} because of the limitations for reproducing important structural three-dimensionality, impacting on cell behavior through adhesive, mechanical, and chemical cues \cite{edmondson2014three}.

Mathematical models and computational tools have also emerged as powerful frameworks for the study of tumor progression, both as a complement of pre-clinical research and as predictive tools in computational oncology \cite{sanga2007predictive,deisboeck2011multiscale,yankeelov2015toward}. Parametric models, capable to assimilate biological data at different scales, are invaluable. The payoff of the mathematical abstractions (after sufficient validation) is to get information on the mechanisms driving tumor progression \cite{altrock2015mathematics},  build predictive tools \cite{enderling2019integrating}, include treatment outcomes \cite{mckenna2018precision}, explore ``what if'' situations, adapt the framework to complex geometries and patient-specific contexts \cite{baldock2013patient}, discard, validate or refine models \cite{hiremath2025identifiability,hormuth2017mechanically}, design new experiments \cite{perez2023massive}, support clinical decision making \cite{perez2024patient}, or integrate multiple scales of analysis \cite{anderson2008integrative,yankeelov2016multi,kazerouni2020integrating}. In fact, mathematical oncology is a valuable tool both when theory guides experiments and experiments guide new theory formulations \cite{byrne2006modelling} and nowadays is an opportunity for personalized medicine \cite{yankeelov2024designing}.

The use of Partial Differential Equations (PDEs) incorporating the specific behavior of the cells of the
tumor under study has proven to be useful to reproduce tumor progression in
experimental heterogeneous cell cultures \cite{ayensa2023mathematical} and in patient-specific
geometries \cite{lorenzo2016tissue}, both of which can be monitored over time while maintaining a
moderate computational cost. Once calibrated, the models can be used to build digital twins \cite{wu2022integrating}
both to predict patient outcomes or to test drugs and treatment schedules \emph{in silico}. This strategy was coined some years ago as \emph{predictive medicine} \cite{oden2016toward}.

Despite the benefits of parametric models, the calibration process is challenging. First,
data quantity and quality are fundamental for parameter identifiability \cite{raue2014comparison}
and calibration \cite{read2020strategies}. Second, complex nonlinear models,
mechanistically representing biological processes are usually multiparametric, thus being prone to
overfitting. In any case, their validation is key, since these models will eventually be used in predictive scenarios \cite{metzcar2024review}. Fortunately, the improvements in cell culture techniques, the development of new biomarkers, and the advent of new and more effective microscopy and imaging techniques triggers the assimilation of much richer although unstructured data \cite{ayuso2022role}.

Recently, enormous efforts have been made to establish a robust methodology for parametric adjustment in biology, particularly computational oncology. The seminal work \cite{oden2016toward} was the starting point for rigorous methodologies culminating in the OPAL method \cite{lima2016selection,oden2017predictive}, which has been successfully applied to gliomas \cite{hormuth2017mechanically} or breast cancer \cite{phillips2023assessing}. In addition, some context-dependent specific algorithms have been developed for learning purposes in real data scenarios \cite{agosti2020learning} and nonconventional approaches have been proposed recently, such as the use of Physics-Informed Neural Networks \cite{daneker2023systems} or hybrid approaches that combine PDEs with some elements of Machine Learning \cite{camacho2025physics}.

However, even if model selection is easily incorporated in most of the aforementioned workflows,
they require both calibration and validation scenarios to decide on the adequacy of a
particular model to be adopted. The concept of model discrepancy or model inadequacy is
underutilized. The existence of a discrepancy term is discussed in many works (see, for
instance, \cite{collis2017bayesian} for a Bayesian approach), but it is not explicitly included and
therefore it is not exploited for decision-making. In addition, the computational cost of PDE models may require ancillary methods to accelerate their simulation for applications that require a high number of model runs, such as uncertainty quantification or reliability analysis \cite{lorenzo2024patient}, where the full \emph{a posteriori} distribution is needed. Therefore, the use of inadequacy-aware surrogates is key. Indeed, in many contexts, we do not need the model to be right, but not to be wrong \cite{enderling2021all}, so the model adequacy has to be evaluated, without having validation scenarios. Our main goal is to establish a framework for biological problems, where parametric calibration can be applied systematically, the model
inadequacy can be evaluated, and whose evaluation might be replaced with that of a surrogate, if needed.

With the idea of exploiting the advantages of Bayesian approaches for inverse problems
\cite{latz2023bayesian}, we present a methodology for model calibration in
biology and biomedical engineering based on the Kennedy and O'Hagan error split
\cite{kennedy2001uh}. The presented approach leverages this error decomposition to evaluate model
adequacy, to assess the use of model surrogates, and to integrate or reconstruct the
probability distribution of the experimental error. In addition, the posterior probability distribution of the parameters does not assimilate neither data noise nor model discrepancy, thus offering an improved biological characterization and therefore impacting downstream predictions, in terms of bias and uncertainty.

We illustrate the approach with a parametric model of spatial \emph{in vitro} cell
progression using real experimental data. This example is relevant enough from the
biomedical point of view: we calibrate a parametric model based on PDEs that attempts to recreate
the progression of glioblastoma (GBM), the most common and aggressive brain cancer
\cite{brat2012glioblastoma,tan2020management} in microfluidic devices, the most biomimetic cell
culture technique \cite{sackmann2014present}. Microfluidic platforms are
very well-suited for recreating GBM physiological and pathological situations for research purposes \cite{ayuso2016development,ayuso2017glioblastoma} and for drug testing \cite{bayona2024tetralol}. In addition, we focus on the parameters related to the go-or-grow metabolic switch \cite{hatzikirou2012go}, a hallmark of GBM progression, using a previous parametric model \cite{ayensa2020mathematical} whose use has been demonstrated to be relevant for incorporating spatial cell cultures techniques in personalized medicine \cite{perez2021predicting}.

This paper is organized as follows. First, in Section~\ref{sec-cell}, we introduce the mathematical model describing glioblastoma cell evolution, which accounts for growth, migration saturation, and oxygen moderation, leading to our case study based on a parametric model of tumor progression. We then present the model calibration framework and discuss the different approaches considered in Section~\ref{sec-problem}. The results in  Section~\ref{sec:res} compare outcomes across four approaches, distinguishing between the use of the parametric versus surrogate model for calibration, and the presence or absence of a model discrepancy term. We also discuss the implications and significance of these findings as a foundation for future work. Finally, we conclude with a summary of our main contributions and propose avenues for further research in Section~\ref{sec:conc}.
\section{Mathematical model of glioblastoma cell culture evolution}
\label{sec-cell}
\begin{figure}[htb]
    \begin{subfigure}[b]{.38\textwidth}
         \centering
         \includegraphics[width=\linewidth]{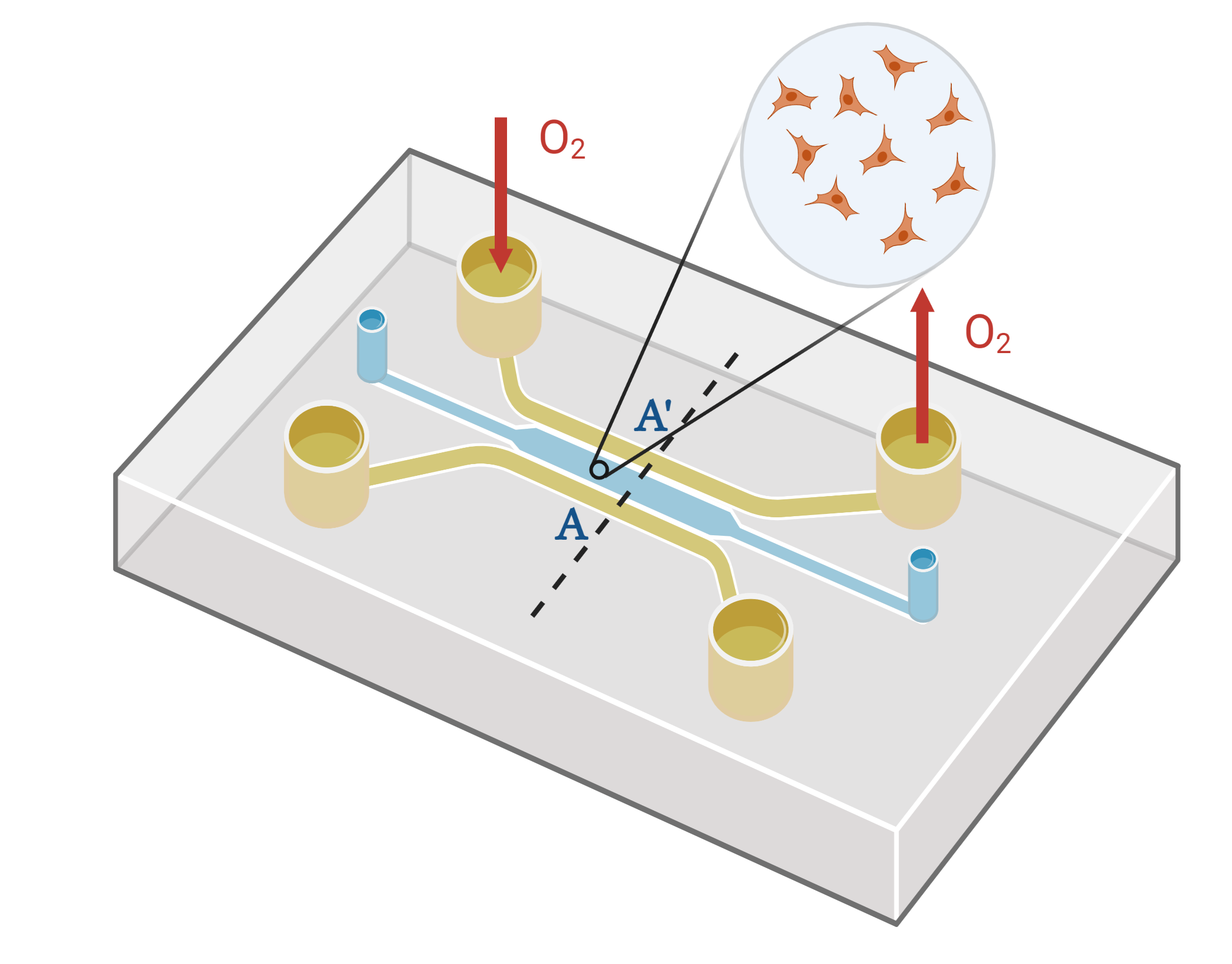}
         \caption{Experimental set-up}
     \end{subfigure}
    \begin{subfigure}[b]{.61\textwidth}
         \centering
         \includegraphics[trim = 0pt 150pt 0pt 100pt, width=\linewidth]{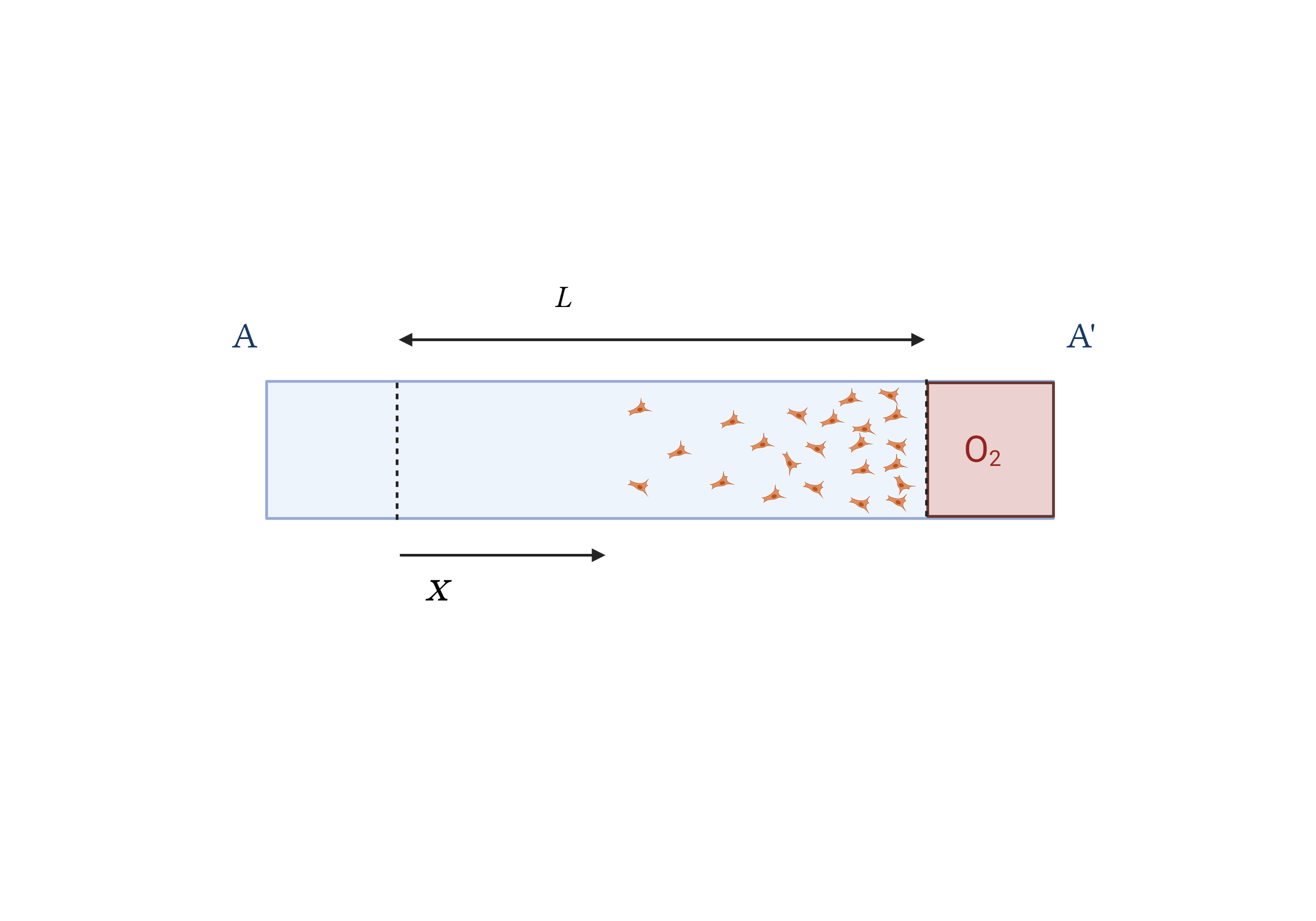}
         \caption{Detail of A-A' section.}
     \end{subfigure}
    \caption{\textbf{Experimental set-up.} A scheme of the microfluidic device is illustrated, together with the geometry and main variables of the model.}
    \label{fig:exp_setup}
\end{figure}

In this section we describe the mathematical model of GBM progression in microfluidic devices. The setup of the experiment is shown in Fig. \ref{fig:exp_setup}.
\subsection{Governing equations}
Our starting point is a nonlinear reaction-diffusion system of PDEs. In the devices under consideration, live and dead cell populations coexist. The PDEs that model their evolution are formulated in terms of the density of live and dead cells ($v$ and $u$ respectively) and of oxygen concentration ($w$). Following the work \cite{ayensa2020mathematical}, the governing equations are:
\begin{subequations} \label{eq::governing}
\begin{align}
    \frac{\partial u}{\partial t} &=  \frac{\partial}{\partial x}\left(D_n\frac{\partial u}{\partial x} - \chi F_\mathrm{go}(u)\Pi_\mathrm{go}(w)u\frac{\partial w}{\partial x}\right) + \frac{1}{\tau_n}F_\mathrm{gr}(u)\Pi_\mathrm{gr}(w)u - \frac{1}{\tau_d}\Pi_\mathrm{d}(w)u, \label{eq::governing_a} \\
    \frac{\partial v}{\partial t} &=  \frac{1}{\tau_d}\Pi_\mathrm{d}(w)u, \label{eq::governing_b} \\
    \frac{\partial w}{\partial t} &=  \frac{\partial}{\partial x}\left(D_{O_2}\frac{\partial w}{\partial x}\right) - \alpha\Pi_\mathrm{c}(w)u, \label{eq::governing_c}
\end{align}
\end{subequations}
where the three unknown fields $u,v,w$ are functions of time~$t\ge0$ and the space coordinate~$x\in[0,L]$.

The first term of the right-hand side of Eq.~(\ref{eq::governing_a}) represents the flow associated
with cell culture migration and has two contributions: the non-oriented motility term
$D_n\frac{\partial u}{\partial x}$ (modeled here as a random diffusion process) and the chemotaxis
term $- \chi F_\mathrm{go}\Pi_\mathrm{go}u\frac{\partial w}{\partial x}$, where cell motility is
induced by the oxygen gradient. $F_\mathrm{go}$ and $\Pi_\mathrm{go}$ are correction factors that
account respectively for the limitations in the motion due to cell crowding and the changes
occurring in cell metabolism associated with migration, depending on the oxygen levels, whose
mathematical structure is discussed later. The second and third term corresponds to the
reaction contribution and it is associated with cell growth and death. The correction terms $F_\mathrm{gr}$,  $\Pi_\mathrm{gr}$ and $\Pi_\mathrm{d}$ have a similar meaning and are also explained later.

With respect to the dead cell equation, Eq. (\ref{eq::governing_b}), it is assumed that dead cells do not migrate, neither by pedesis nor by chemotaxis. Therefore, we only consider the source term coming from the death of the cells.

The third equation,  Eq. (\ref{eq::governing_c}), refers to the oxygen evolution equation; the first
term of the right-hand side models the flow contribution, consisting solely of oxygen diffusion, and
the second corresponds to the oxygen consumption by cells. As for the cell equation, the function $\Pi_\mathrm{c}$ is a nonlinear correction function.
\subsection{Boundary and initial conditions}
Eqs. (\ref{eq::governing_a}), (\ref{eq::governing_b})  and (\ref{eq::governing_c}) must be
complemented with appropriate boundary and initial conditions, which are the data at a specific time. Note that we
will refer to this time as $t=0$ even if it is not necessarily the experiment starting time, but the instant when the cell culture profile is measured with the microfluidic device
fully oxygenated.
\begin{subequations} \label{eq::initial}
\begin{align}
    u(x,t=0) &= u_0(x), \\
    v(x,t=0) &= v_0(x), \\
    w(x,t=0) &= w_0,
\end{align}
\end{subequations}
with $u_0$ and $v_0$ given known function and $w_0$ the (homogeneous) ambient oxygen level. 

With respect to the boundary conditions, we assume that some cells can escape from the central chamber to the lateral channels, but this cell flow is proportional to the cell concentration at the interface. No migration is considered for dead cells, so no boundary conditions associated with Eq. (\ref{eq::governing_b}) are necessary. Regarding oxygen evolution, the oxygen concentration at the boundaries is assumed to be known and equal to $w_0$, that is, the ambient concentration. Therefore, the boundary conditions for live cell population are:
\begin{subequations} \label{eq::boundary}
\begin{align}
    u(x=0,t) - j\,f_u(x=0,t) &= 0, \\
    u(x=L,t) + j\,f_u(x=L,t) &= 0,
\end{align}
where $j$ is a parameter that controls the proportionality between cell concentration and cell flow given by
$$f_u(x,t) = D_n \frac{\partial u}{\partial x}(x,t) - \chi F_\mathrm{go}(u(x,t))\Pi_\mathrm{go}(w(x,t))u(x,t)\frac{\partial w}{\partial x}(x,t),$$
whereas for oxygen concentration
\begin{align}
    w(x=0,t) &= w_0, \\
    w(x=L,t) &= w_0.
    \end{align}
\end{subequations}
\subsection{Model corrections}
To complete the model, we need to define the corrections that are specific to GBM, that is, the functions $F_\mathrm{go}$, $F_\mathrm{gr}$, $ \Pi_\mathrm{go}$, $\Pi_\mathrm{gr}$, $\Pi_\mathrm{d}$ and $\Pi_\mathrm{c}$
that encode the cell behavior.
\subsubsection{Growth and migration saturation}
The functions $F_\mathrm{gr}$ and $F_\mathrm{go}$ are related to saturation considerations as cell growth and migration are limited by the amount of space. Therefore, $F_\mathrm{go}$ and $F_\mathrm{gr}$ depend on the cell density $u$. Following the work \cite{ayensa2020mathematical}, we write
\begin{subequations}
\begin{equation}
    F_\mathrm{gr}(u,v) = 1 - \frac{u + v}{c_\mathrm{sat}},
\end{equation}
\begin{equation}
    F_\mathrm{go}(u) = 1 - \frac{u}{c_\mathrm{sat}},
\end{equation}
\end{subequations}
where $c_\mathrm{sat}$ is the cell saturation.
With respect to the nonlinear growth correction, it is required because dead cells occupy physical
space, restricting the ability of viable cells to divide and also that cellular debris and necrotic
material released by dead cells may alter the local microenvironment by generating toxic
by-products, disrupting pH homeostasis, and affecting signaling pathways essential for cell growth
\cite{hanahan2022hallmarks}. Chemical cues released by dead cells inhibiting migration are neglected, hence
only space needs are considered which are mainly modulated by live cell density.
\subsubsection{Oxygen moderation}
For some cellular mechanisms, oxygen acts as a moderator variable. In fact, cell hypoxia has been shown to play an important role in cell culture dynamics \cite{brat2004pseudopalisades,rong2006pseudopalisading,ayuso2016development,ayuso2017glioblastoma}, particularly in cell proliferation, death, and migration. These hypoxia-mediated dynamics have been proposed as one of the main driving forces of GBM invasion and aggressiveness~\cite{brat2004vaso}.

\paragraph{The go-or-grow paradigm} The metabolic behavior of the GBM cells, in particular, its response to changes of oxygen pressure, is encoded in the functions $\Pi_\mathrm{go}$ and $\Pi_\mathrm{gr}$. These two functions regulate the activation/deactivation of migration and proliferation. There is sound evidence in the scientific literature that the switch between proliferative and migratory activity in a cell population is hypoxia-mediated \cite{lu2010hypoxia,carreau2011partial}, that is, $\Pi_\mathrm{go} = \Pi_\mathrm{go}(w)$ and $\Pi_\mathrm{gr} = \Pi_\mathrm{gr}(w)$. We know that $\Pi_\mathrm{gr}:\mathbb{R}^+ \rightarrow [0;1]$ is a continuous, non-decreasing
function and $\Pi_\mathrm{go}:\mathbb{R}^+ \rightarrow [0;1]$ is another, continuous non-increasing
function. A slightly more restrictive condition is to consider that, in addition to monotonicity, we
must have
$$\lim_{u \rightarrow 0^+} \Pi_\mathrm{go}(w) = \lim_{u \rightarrow +\infty } \Pi_\mathrm{gr}(w) = 1, \qquad  \lim_{u \rightarrow +\infty} \Pi_\mathrm{go}(w) = \lim_{u \rightarrow 0^+ } \Pi_\mathrm{gr}(w) = 0.$$

Finally, based on energetic considerations, we may assume that cells invest their resources either in proliferation or migration, so it is plausible to postulate that there exists a function $Q : [0,1] \times [0,1] \rightarrow \mathbb{R}$ such that $\forall w \in \mathbb{R}^+$,
\begin{equation*}
Q(\Pi_\mathrm{go}(w),\Pi_\mathrm{gr}(w)) = 1.
\end{equation*}

In the work by \cite{ayensa2020mathematical}, a piecewise linear function is proposed to
represent the go-or-grow behavior, namely,
\begin{subequations}\label{eq:go_or_grow}
\begin{equation}\label{eq:grow}
    \Pi_\mathrm{gr}(w) =
    \left\{
    \begin{array}{lcc}
        \frac{w}{h_1} &   \mathrm{if}  & 0 \leq w \leq h_1 \\
        1 &  \mathrm{if} &  w  > h_1
    \end{array}
    \right. ,
\end{equation}
\begin{equation}\label{eq:go}
    \Pi_\mathrm{go}(u) =
    \left\{
    \begin{array}{lcc}
        1-\frac{w}{h_1} &   \mathrm{if}  & 0 \leq w \leq h_1 \\
        0 &  \mathrm{if} &  w  > h_1
    \end{array}
    \right. .\\
\end{equation}
\end{subequations}
Here, $h_1$ is the hypoxia threshold, the oxygen characteristic scale below which the cell metabolism changes from a proliferative to migratory nature. It is clear that for the considered model, the function that takes into account energy competing resources is $Q(x,y) = x + y$.

\paragraph{Cell death: apoptosis and necrosis} Cell death is a natural process that depends on many
factors and agents with intrinsic stochastic nature \cite{galluzzi2018molecular}. Anoxia is one
fundamental cause of cell death \cite{sendoel2014apoptotic}. We can therefore assume that
$\Pi_\mathrm{d} : \mathbb{R}^+ \rightarrow [0,1]$ is a non-increasing function. However, since
apoptosis is always a cell death mechanism, it is reasonable to assume that $\forall w \in
\mathbb{R}^+, \ \Pi_\mathrm{d}(w) > 0$. In addition, we can impose that
$$\lim_{w \rightarrow 0^+} \Pi_\mathrm{d}(w) = 1\ ,
\qquad
\lim_{w \rightarrow +\infty} \Pi_\mathrm{d}(w) = 0\ .
$$
The work \cite{ayensa2020mathematical} proposes the following nonlinear
correction:
\begin{equation}\label{eq:anoxia}
    \Pi_\mathrm{d}(w) = \frac{1}{2}\left[1 - \tanh\left(\frac{w - h_2}{\Delta h_2}\right)\right].
\end{equation}
In Eq. (\ref{eq:anoxia}), $h_2$ is the anoxia threshold and $\Delta h_2$ the anoxia sensitivity. These two parameters regulate the impact of oxygen on the metabolism of cell death, that is, whether the cell dead is due solely to apoptosis or to both apoptosis and necrosis.

\paragraph{Oxygen consumption} Finally, oxygen consumption is related with the kinetics of oxidative phosphorylation that occurs in the membrane of cellular mitochondria, usually described using a non-decreasing function  $\Pi_\mathrm{c} : \mathbb{R}^+ \rightarrow [0,1]$ such that
$$\lim_{w \rightarrow 0^+} \Pi_\mathrm{c}(w) = 0\ ,
\qquad
\lim_{w \rightarrow +\infty} \Pi_\mathrm{c}(w) = 1\ .
$$
It is usual in biology and pharmacology to use a Hill function for this purpose
\cite{goutelle2008hill}. Following again the work 
\cite{ayensa2020mathematical}, we assume Michaelis-Menten kinetics for oxidative phosphorylation, thus
$$\Pi_\mathrm{c}(w) = \frac{w}{w + k_m}.$$
\subsection{Case study: parametric model of glioblastoma progression}\label{subs-param_model}
At this point, the presented model depends on eleven parameters, namely, $D_n$, $D_{O_2}$, $\chi$,
$\tau_n$, $\tau_d$, $\alpha$, $c_\mathrm{sat}$, $h_1$, $h_2$, $\Delta h_2$, $k_\mathrm{m}$. Some of
them have a well-identified value in the scientific literature. For example, for those associated to
the oxygen evolution equation, the oxygen diffusion $D_{O_2} = 1 \times 10^{-5} \, \mathrm{cm^2\cdot s^{-1}}$ has been reported in many works \cite{tannock1972oxygen,dacsu2003theoretical} as well as the Michaelis-Menten constant, $k_\mathrm{m}=2.5 \, \mathrm{mmHg}$ \cite{alper1956role,dacsu2003theoretical}. For some others, the value is not so well characterized, although simple conventional experiments could be used for calibration purposes:
\begin{itemize}
    \item The parameters related to cell growth and death under fully oxygenated conditions, $\tau_n$, $\tau_d$ and $c_\mathrm{sat}$, can be determined in cell growth experiments in ambient conditions and in the absence of oxygen gradients, both in microfluidic devices \cite{lei2014real,tao2015based}, using cell spheroids \cite{vinci2012advances}, or even in conventional Petri dish experiments.
    \item Assuming that $k_m$ is known, the oxygen consumption rate, $\alpha$ is easily obtained by measuring the ambient oxygen pressure in an isolated system with a controlled cell culture population. It is even possible to determine both $k_\mathrm{m}$ and $\alpha$ from the oxygen pressure using an Eadie–Hofstee diagram \cite{hofstee1959non} or a Lineweaver–Burk plot \cite{lineweaver1934determination}.
    \item The identification and calibration of the parameters related to the metabolic changes of cell death mechanism, that is, $h_2$ and $\Delta h_2$ is more subtle, since it requires the experimental measurement of both the ambient oxygen pressure and the dead cell density evolution during the experiment.
    \item Finally, the non-oxygen-mediated pedesis constant, $D_n$, is more complicated to determine as spatial cell cultures are necessary. However, spheroid cultures \cite{ayuso2015study} and microfluidic devices \cite{mehling2014microfluidic} offer a great opportunity for cell migration evaluation. If full oxygenation is guaranteed in the whole culture, no oxygen gradient is formed so non-oxygen mediated pedesis is easily computed, for instance, once given $\tau_n$, $D_n$ can be determined by evaluating the cell migration radial velocity $V$ and using the Fisher's model \cite{fisher1937wave}, $V = 2\sqrt{D_n \tau_n}$.
\end{itemize}
The values of $\chi$ and $h_1$ are nevertheless significantly more difficult to determine. To
identify their values we need to experimentally establish heterogeneous oxygen distributions. Indeed, as the cell migration depends on the oxygen level (and not only on the oxygen gradient), it is difficult to estimate this value from a single experiment or measurement. However, we can measure cell culture migration under an oxygen gradient in a very localized region where the oxygen level may be considered almost constant \cite{funamoto2012novel}. From now on, we will define $b = 1/h_1$ , the inverse of the hypoxic threshold, as the model parameter for comparison purposes.

In addition to the model parameters, we have two constants $j$ and $w_0$ associated with the boundary
conditions, assuming the initial profiles $u_0(x)$ and $v_0(x)$ as given. The ambient oxygen
condition $w_0$ can be easily measured, but the parameter $j$ is related to cell fluxes and depends both on the properties of the cell culture gel, the physics of the interface lateral/central chambers, and the cell properties, so it is difficult to fix it \emph{a priori}.

Therefore, we consider the following parametric model $\mathcal{M}_{\bs{\theta}}$:
\begin{equation}
    \eta(x;\bs{\theta}) = \hat{u}(x,t=T;\bs{\theta}),
\end{equation}
where $\hat{u}$ is the live cell density field obtained by solving numerically the PDEs (\ref{eq::governing}) with initial conditions (\ref{eq::initial}) and boundary conditions (\ref{eq::boundary}),  $T$ is a given experimental and simulation horizon and the (unknown) parameter vector is $\bs{\theta} = (\tau_n,\chi,b,j)$. The reason for including the parameter $\tau_n$ is because the triad of parameters $(\tau_n,\chi,h)$ fully characterize the \emph{go-or-grow} metabolic switch \cite{hatzikirou2012go}, that has been the subject of scientific research (see for instance \cite{ayensa2025overview} for an extensive review). The remaining model parameters are the ones that are considered in \cite{ayensa2020mathematical}.
\section{Model calibration}
\label{sec-problem}
The model described in Section~\ref{sec-cell} can be employed to predict cell migration, provided
that the boundary and initial conditions are set up, and all the parameters appearing in the
equations are given appropriate values. This last aspect, however, entails many more difficulties
that might naively appear. Moreover, an incorrect parameter selection can compromise the modeling and
computational efforts spent towards migration predictions.
\subsection{Calibration approaches}
In general, determining the correct values for the parameters of a model is not trivial. These must be adjusted using
information obtained from experimental data but, since the latter do not provide direct
parameter values, some type of \emph{inference} must be performed. We refer to as \emph{calibration}
to the process of selecting the best values for the parameters given the information available.
Naturally, since there is not a single choice for what is \emph{best} and different sources of
information can be incorporated or disregarded, several alternative and valid calibration strategies exist.

The simplest and most common approach to model calibration is based on mismatch minimization: given
some norm (for example, the Euclidean norm) measuring the discrepancy between the model predictions
and the measurements on a set of experiments, the parameters are selected to minimize this error
(see, e.g., \cite{hansen2014ye} for a comprehensive overview). Even if often employed, this
calibration strategy has well-known limitations: it is very sensitive to outliers, it might require
some regularization that affects the solution, it might select parameter values that only yield local
minima to the mismatch error, and, finally, can only provide point estimates possibly with often very wide confidence intervals for the optimal parameter values~\cite{backman2023}.

More powerful methods for calibration can be obtained by leveraging Bayesian statistics. Building on
the classical work of \cite{kennedy2001uh}, numerous applications have appeared
in the literature that refine the original model \cite{higdon2004dl,van2005kg,ohagan2006cm,higdon2008in} or apply it to various fields
\cite{chong2018pj,willmann2022jb,depablos2023ub}. The Bayesian point of view corrects many of the
deficiencies of the least-squares approach, being especially suitable when the number of experiments
that can be conducted is scarse.

There is a fundamental difference between the minimization approach to calibration and the Bayesian
one. In the first class of methods, optimal \emph{values} are sought for the unknown parameters. In
the latter, optimal \emph{probability distributions} are the object of interest. In this statistical
picture, it is therefore accepted from the outset that parameters of (almost) all models are random variables. The goal of Bayesian calibration is not to find a single value for
each parameter, but rather understand the probability distribution of each value, the correlations between pairs of parameters, etc. This probabilistic point of view
is indispensable to tackle all the critical issues associated with the \emph{uncertainty} of models
and their predictions.

In this work, and in contrast with our previous article \cite{ayensa2020mathematical}, we embrace
the Bayesian approach to calibration of models for cell migration, and provide a detailed
exploration of its benefits. We start next by summarizing its main features in an abstract fashion before applying it to the models of interest in this work.
\subsection{Bayesian calibration}
\label{subs-bayesian}
To start, let us consider the output of a physical phenomenon (or economical, chemical, social, etc.) which is denoted as $\varphi$ and is a function of some controllable variables, an array $\bs{x}\in\mathbb{R}^k, k\ge1$. When a measurement~$z$ is made for $\varphi(\bs{x})$, the result obtained is never exact and thus we can write
\begin{equation}
  \label{eq-experimental}
  z = \varphi(\bs{x}) + \epsilon(\bs{x})\ ,
\end{equation}
where $\epsilon$ is a random variable accounting for the experimental error, usually assumed to be normally distributed. From an experimental campaign, one often has access to a dataset $\mathcal{D}=\{(\boldsymbol{x}_i,z_i)\}_{i=1}^M$ of $M$ (inexact) measurements. 

Let us now consider an analytical or computational model $\eta:\mathbb{R}^k\times\mathbb{R}^p\to
\mathbb{R}$ of the same phenomenon as above where $\eta(\bs{x};\bs{\theta})$ is now a function of the
controllable variables~$\bs{x}$ and an array of parameters~$\bs{\theta}\in \mathbb{R}^p$. In contrast with the previous situation, now the model $\eta$ is completely known; in turn, the parameters are either unknown or, more commonly, known only approximately. This situation is the one we face with the cell migration model of Section~\ref{sec-cell}.

The Bayesian approach to model calibration starts from the assumption that searching directly for the optimal value of the model parameters~$\theta$ is inconsistent with our knowledge and the information provided by the experiments. Rather, based on the prior information that one might have about these parameters and the experimental results, there exists a probability distribution for $\bs{\theta}$ that best explains the measurements, incorporating the previously existing information. The ultimate goal of \emph{Bayesian inference} is the determination of the latter probability distribution, the so-called \emph{posterior} distribution.
\subsection{The Kennedy and O'Hagan split}
Next, we adapt the Kennedy and O'Hagan (KOH) ideas for calibrating a computer code that solves an
initial boundary value problem of the type describing cell migration. The starting point of this
Bayesian strategy is, as anticipated, an \emph{a priori} multivariate probability distribution for the parameters whose probability density function will be denoted as $\pi=\pi(\bs{\theta})$. This distribution should reflect, as best as possible, the previous knowledge that we have about the model parameters, be it from previous models, other research groups, preliminary tests,~etc.

To \emph{infer} a better probability distribution for the parameters, any Bayesian method requires two more ingredients: data and a model for the likelihood of these data as a function of the parameters. The first one is provided by an experimental campaign $\mathcal{D}$ as defined above. To construct the second and more complex element let us note, first, that the model $\eta$ is typically insufficient to predict correctly all the experimental results and thus we assume that
\begin{equation}
  \label{eq-koh}
   z_{i} = \eta(\boldsymbol{x}_{i}; \bs{\theta}) + \delta(\boldsymbol{x}_{i}) + \epsilon(\boldsymbol{x}_{i})\ ,\qquad i=1,2,\ldots,M\ ,
\end{equation}
where $\delta:\mathbb{R}^k\to \mathbb{R}$ is the \emph{discrepancy} error, an unknown function that accounts for the inability of the model $\eta$ to reproduce $\varphi$, even with the best choice of parameters. The objective of the KOH calibration is thus not only to obtain the best \emph{a posteriori} distribution for $\bs{\theta}$, but also to reveal~$\delta$, at least in a way that is as consistent with the data and prior information as possible.

A convenient way to work with a general class of \emph{non-parametric} functions and select from them the discrepancy error is to assume it is a Gaussian process. In addition, the model $\eta$ might be expensive to evaluate so it also proves convenient to introduce a surrogate or meta-model~$\tilde{\eta}$ for it, one that can be evaluated much faster. Since the error is normally distributed and the discrepancy is a Gaussian process, the calculations become simple when also $\tilde{\eta}$ is selected among Gaussian processes.
\subsection{Gaussian processes}
Gaussian processes (GP) play a key role in the calibration framework proposed by Kennedy and O'Hagan. GPs can be understood as generalizations of Gaussian multivariate distributions. Specifically, they are a collection of random variables indexed by a real number such that every finite subset of these variables has a multivariate normal distribution. We indicate that $f$ is a GP with mean function $m$ and covariance $c$ using the notation $f\sim \mathcal{GP}(m(\bs{x}),c(\bs{x},\bs{x}'))$. See \cite{rasmussen2006vz} for a comprehensive description of these concepts.

In practice, here, we will always employ GPs with $m\equiv 0$ and covariance functions of the form
\begin{equation}
    c(\bs{x},\bs{x}') = \lambda \exp\left[-\frac{\|\bs{x}-\bs{x}'\|^2}{2\beta^2}\right] ,
\end{equation}
with $\lambda, \beta$ referred to as \emph{hyperparameters}, $\lambda$ denoting the signal variance and $\beta$ the lengthscale. This corresponds to a stationary isotropic squared exponential kernel.
In the calibration, the hyperparameters are unknown and their optimal probability distributions will
have to be obtained, as for any other unknown parameter of the model.

In particular, the hyperparameter $\beta$ has the same dimensions as the input variable $\bs{x}$, as it
is thus referred to as a \emph{lengthscale}. Note, also, that the norm of the input variable $\bs{x}$ is only
well-defined if the dimensions of all its components are the same. It therefore becomes practical to
nondimensionalize the model equations from the outset, ensuring the consistency of the norm operations.
\subsection{Calibration procedure}
The Bayesian calibration consists, essentially, in calculating a posterior probability distribution given the prior probability distribution of the parameters and the hyperparameters, plus the likelihood function of the data.

In this work, we consider four different calibration procedures that illustrate the main features of Bayesian calibration. In the following, we succinctly introduce them and refer to the literature for more details on the actual calculations involved in their computation \cite{kennedy2001uh,ohagan2006cm,higdon2008in}.

\paragraph{Bayesian inference (BI)} The simplest calibration strategy is a prototypical example of Bayesian
inference. In this case, the discrepancy of the model is ignored, but we have some prior knowledge for the parameters of the model and the measurement error in experimental data, which is assumed to be normally distributed. Then, a straightforward application of Bayes' theorem provides the posterior probability distribution for all the unknown random variables, that is, the parameters and the experimental error. 

Although some analytical results enable the calculation of the posterior for simple priors and low dimensions, for the
sake of generality, this probability distribution is often approximated using MCMC methods (see, for example, general monographs on the subject \cite{robert2004gp,barbu2020is}, or a recent work with references to some of the latest algorithms \cite{romero2025qk}). Let us note that these numerical techniques are extremely efficient, but demand a large number of likelihood and prior evaluations.

\paragraph{Bayesian calibration of expensive models (BCE)}\label{BCE}
It might turn out that the model to be calibrated is computationally expensive and thus impossible
--- or impractical --- to run thousands of MCMC steps that explore the parameter and hyperparameter
space. In these situations, the user might decide to replace the model with a \emph{surrogate} or
\emph{meta-model} that, once trained, can be evaluated with a small cost and used
for calibration. Since the true model is replaced with a surrogate, naturally, the calibration will be less accurate than the one described before.

While replacing an expensive model with a fast surrogate seems appealing, one should be aware of its
limitations. To construct the surrogate, a small set of model evaluations needs to be used in the
training phase. For this set to be representative of the model, the input and parameter spaces need
to be sufficiently covered. In high-dimensional spaces, this methodology is very costly and efficient
sampling methods such as the Latin hypercube are employed \cite{mckay1979xs,iman1981nxy}.
Note that although Latin hypercube sampling is specifically designed to avoid clustering and ensure
a uniform distribution of samples across each dimension, it may suffer from clustering effects when the dimensionality is high. This issue becomes even more apparent in the presence of constraints \cite{Esposito2023,schenk2024castroefficientconstrained}.
Moreover, these methods cannot compensate for a vague definition of the parameter range, especially if
the model is very sensitive to some of them. When the input and parameters are poorly sampled, we
cannot expect the surrogate to be accurate, and the probability distributions of the model parameters
will be uninformative. In these situations, it often happens that the posterior parameter
distributions are almost identical to the prior distributions provided in the calibration, and the
MCMC acceptance ratio is small. These two are signs that the likelihood is very small in the MCMC
chains, and that the posterior is controlled by the prior distributions, meaning that the provided data is uninformative. If any of these flags is raised,
the analyst is advised to narrow the parameter range in the surrogate training, if possible, or
abandon this strategy altogether.

To calibrate a model using a surrogate, one needs to first generate a set of \emph{synthetic} data
points, obtained by evaluating the model with input and parameter values
$\{\tilde{\boldsymbol{x}}_i,\tilde{\boldsymbol{\theta}}_i\}_{i=1}^M$ randomly selected in input and parameter space using, as
indicated, e.g., the Latin hypercube method.

\paragraph{Bayesian calibration with discrepancy (BCD)}\label{BCD}
When calibrating models, it is not uncommon that the optimal parameters result in a model whose
predictions are not sufficiently good. This could be due to a bad model choice, often the price for
selecting a model that is too simple. One of the outcomes of this situation is that the surrogate
model variant (cf. previous subsection) can represent the experimental data even better than the
physical model itself. Even then, the discrepancy of the model can still be estimated. Calibrations of this type do not typically improve the calibrated model; however, they can quantify the mismatch of the model and the predictions, hinting at the need to select a better model if more accuracy is desired.

\paragraph{Bayesian calibration of expensive models with discrepancy (BCED)}
When evaluating a computationally expensive model, and the goal is to assess the quality of the model, here represented by a surrogate, both the physical model and the discrepancy are approximated using Gaussian Process (GP). This allows for quantifying the mismatch between the surrogate and the resulting predictions, thereby informing the selection of more effective synthetic data to construct such a surrogate and/or a more appropriate kernel representation. Moreover, incorporating the discrepancy explicitly can increase the predictive uncertainty, which transparently reflects the combined uncertainties from both the surrogate approximation and possible model inadequacies, ultimately leading to more cautious and reliable predictive statements.
\section{Results and discussion}
\label{sec:res}
Next, we present the main findings from the four Bayesian calibration strategies
described in Section~\ref{sec-problem} evaluating their performance and discussing how the probabilistic results can be leveraged for downstream analyses and decision-making.  Specific data and model configurations can be found in the Appendix. 
\subsection{Evaluating the calibration procedures}
\begin{figure}[htb]
    \centering
    \includegraphics[width=.65\linewidth]{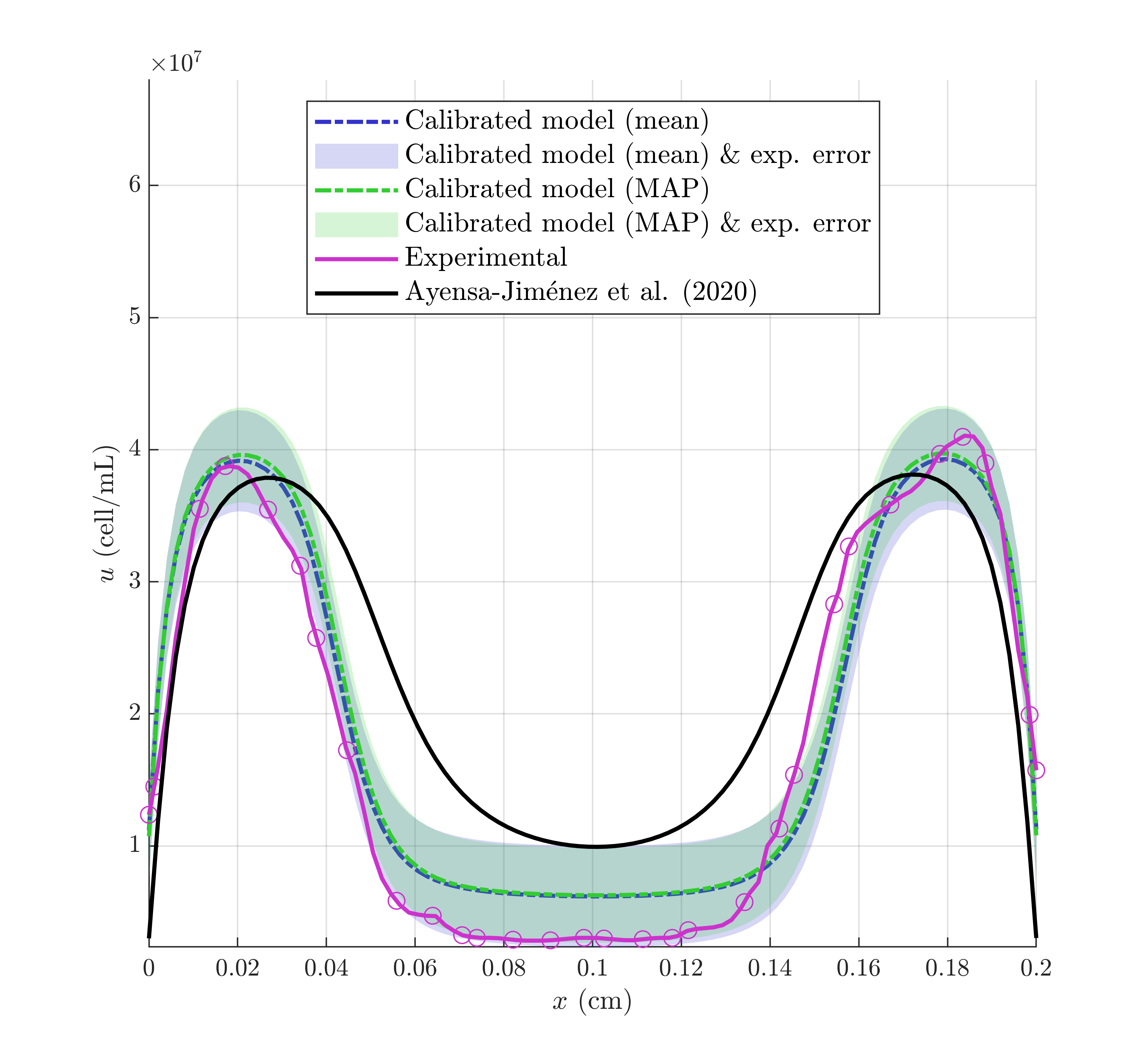}
    \caption{\textbf{Live cell profile generated. BI approach.} Earlier results  \cite{ayensa2020mathematical} showed an error of $\epsilon = 0.0768$, when compared to the prediction using the MAP parameters, with $\epsilon = 0.0326$.}
    \label{fig:M3_TA}
\end{figure}

\begin{figure}[htb]
    \centering
    \includegraphics[width=.65\linewidth]{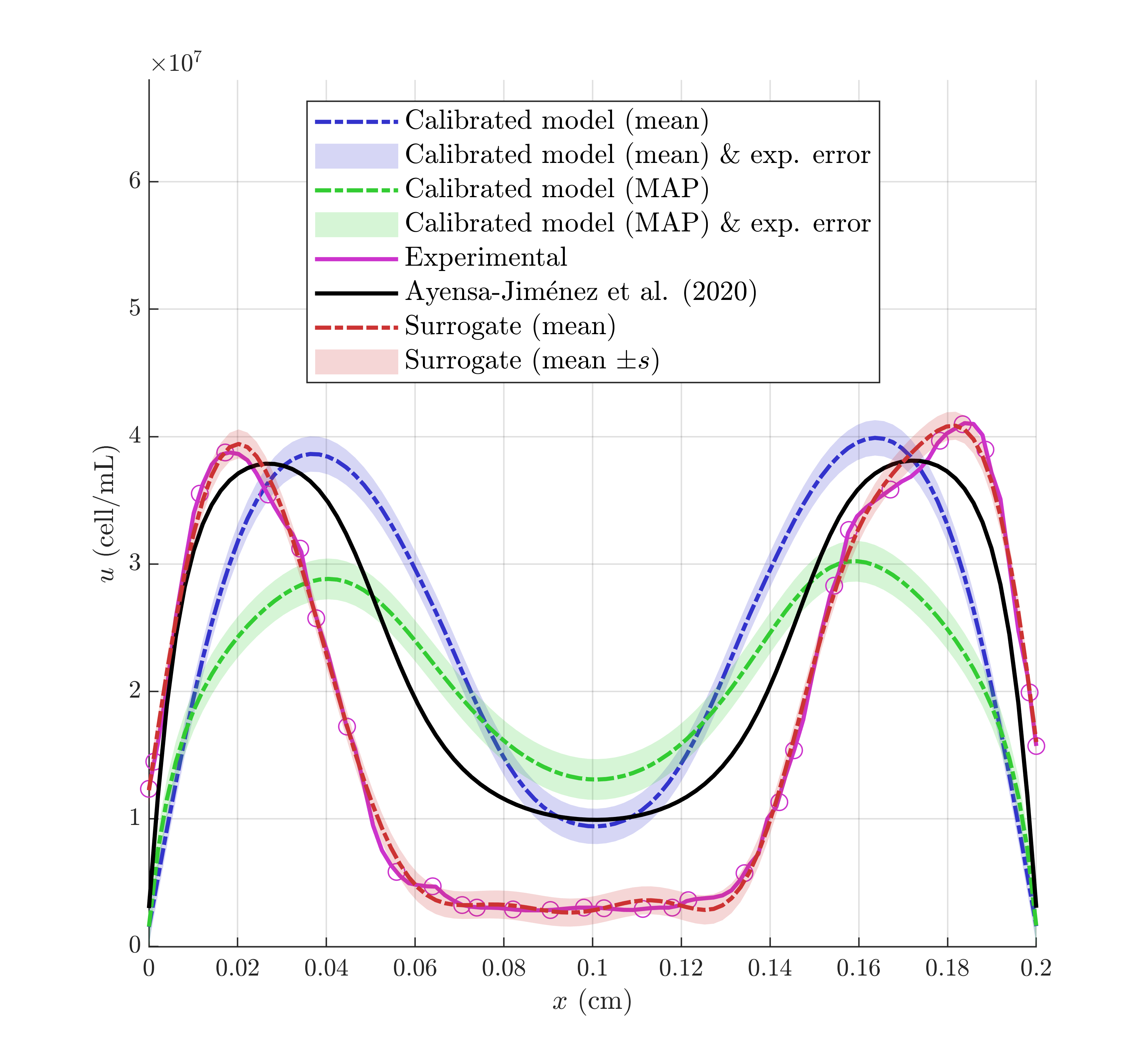}
    \caption{\textbf{Live cell profile generated using BCE approach.} Earlier results  \cite{ayensa2020mathematical} showed an error of $e = 0.0768$, when compared to the prediction using the MAP parameters, with $e = 0.1161$ and the prediction of the surrogate $e = 0.0075$.}
    \label{fig:M3_TB}
\end{figure}

\begin{figure}[htb]
    \centering
    \includegraphics[width=.65\linewidth]{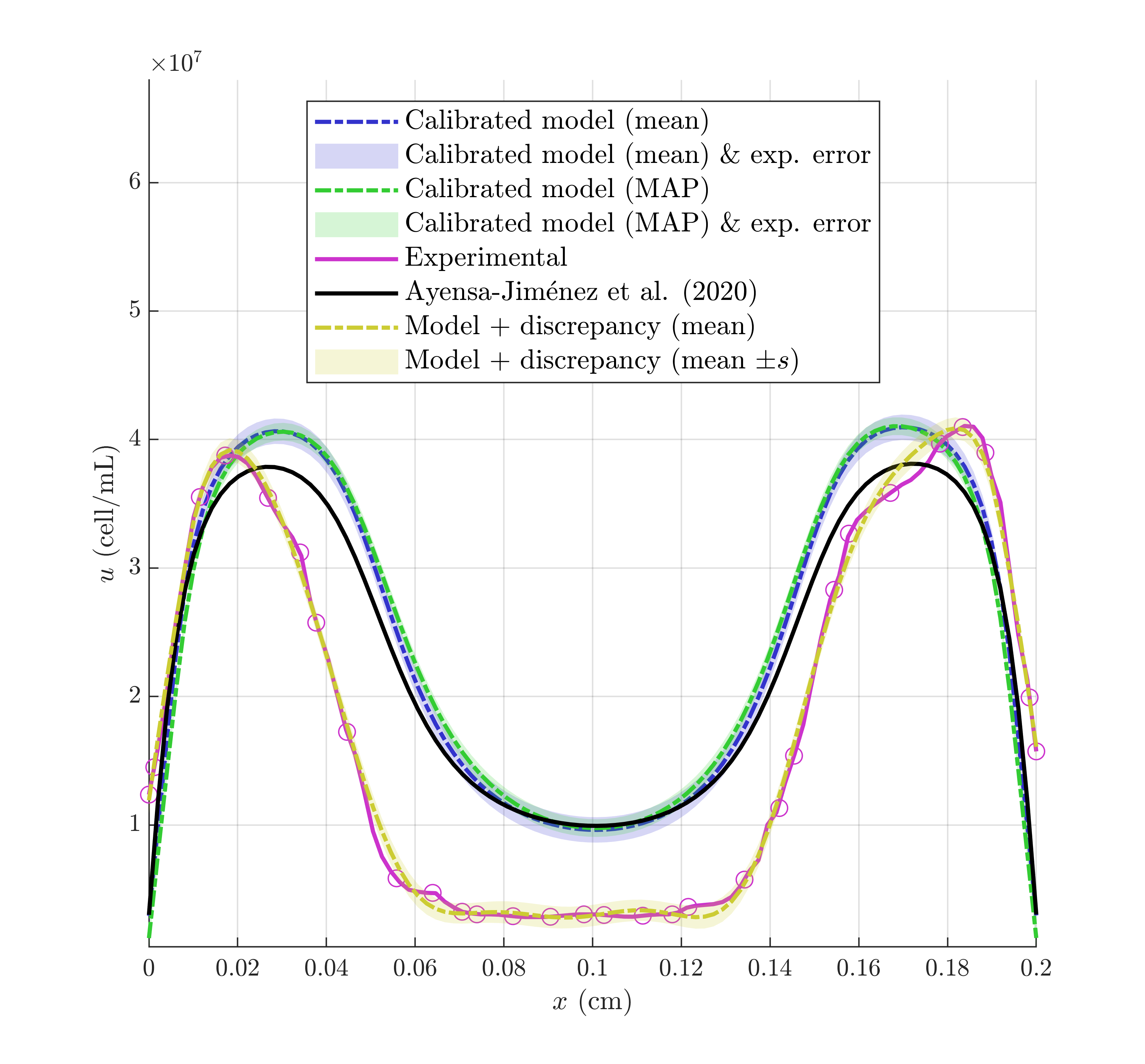}
    \caption{\textbf{Live cell profile generated using the BCD approach.} Earlier results  \cite{ayensa2020mathematical} showed an error of $e = 0.0768$, when compared to the prediction using the MAP parameters, with $e = 0.0948$ and the prediction including the discrepancy $e = 0.0068$.}
    \label{fig:M3_TC}
\end{figure}

\begin{figure}[htb]
    \centering
    \includegraphics[width=.65\linewidth]{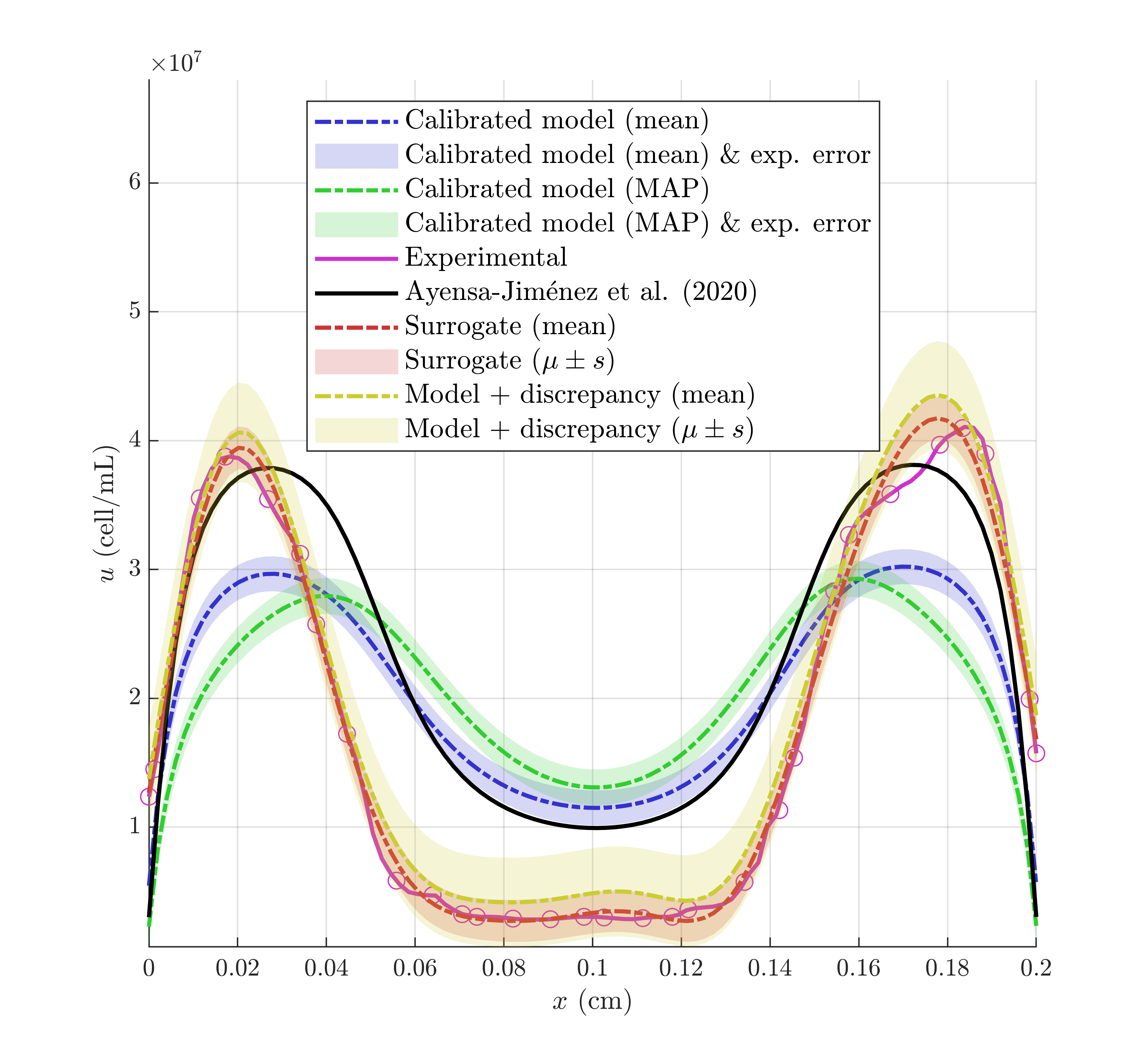}
    \caption{\textbf{Live cell profile generated using the BCED approach.} Earlier results  \cite{ayensa2020mathematical} showed an error of $e = 0.0768$, when compared to the prediction using the MAP parameters, with $e = 0.1132$, the prediction including the discrepancy with $e = 0.0138$ and the one using the surrogate with $e = 0.0116$.}
    \label{fig:M3_TD}
\end{figure}

We first evaluate the calibration methodology using the parametric model proposed in \cite{ayensa2020mathematical} (see Section \ref{subs-param_model}). The results are shown in Fig. \ref{fig:M3_TA} for the BI approach, Fig. \ref{fig:M3_TB} for the BCE approach, Fig. \ref{fig:M3_TC} for the BCD approach and Fig. \ref{fig:M3_TD} for the BCED approach and the corresponding posterior distributions in Fig. \ref{fig:corner_typeA}, \ref{fig:corner_typeB}, \ref{fig:corner_typeC}, \ref{fig:corner_typeD}, respectively. For comparison purposes, we define the prediction error of the live cell profile as
\begin{equation}
    e = \frac{1}{c_\mathrm{sat}}
    \left(\int_0^L \left(\hat{u}(x) - u_\mathrm{exp}(x)\right)^2\, \mathrm{d}x\right)^{1/2}.
\end{equation}
In Fig. \ref{fig:M3_TA} we show the live cell profile obtained when the model parameters and the
experimental error are adjusted using the BI approach and the experimental data. We show
the predicted cell profiles and the experimental error by both MAP (Maximum a posteriori) and (sample) mean estimation of the parameters from the posterior distributions. The predictions obtained using the MAP give an error of $e = 0.0326$, improving the prediction in \cite{ayensa2020mathematical}, with $e = 0.0768$. The mismatch between experimental data and model predictions is mostly absorbed by the experimental error, which is not negligible ($\sigma = 7.2 \times 10^{-2} c_\mathrm{sat}$, so $7.2\%$ of the saturation concentration). Interestingly, the predictions using the MAP and mean values of the parameters are very similar, illustrating the low manageable variability of the posterior parameter distributions.

In Fig. \ref{fig:M3_TB} we show the results obtained when using the BCE approach. In addition to the predictions obtained with the full model, we include the mean $\pm$ standard deviation of the surrogate. Although the surrogate reduces the error to $e = 0.0075$, the error of the predictions obtained using the MAP is $e = 0.1161$, illustrating the failure to match the parametric model and the experimental data in the surrogate. Unlike in the previous case, the experimental error estimation is drastically reduced, at the cost of greatly increasing the variability of the model parameters as suggested by the differences when estimating the cell profile using the MAP and the mean of the parameters.

In Fig. \ref{fig:M3_TC}, we show the predictions of the cell profile using the BCD approach. With the inclusion of the discrepancy, the predictions are of the same order as
ones reported in previous works, with $e = 0.0948$ achieved using the MAP value of the parameters.
The error is reduced to $e = 0.0068$ when considering the model discrepancy. The difference between these latter two errors is an indication of model inadequacy. However, when compared to the previous approaches, both the experimental error and parameter variability (as indicated by the overlap of MAP and mean predictions) are kept small. The mismatch between experimental data and model predictions is fundamentally absorbed by the discrepancy term. This is, indeed, the main component of the model that is able to detect model inadequacy, regardless of experimental error, while keeping the parameter approach valid enough, that is, assuming a posterior distribution of the parameters informative enough.

Finally, in Fig. \ref{fig:M3_TD} we show the predictions using the BCED approach. The error obtained
using the MAP values is $e = 0.1132$, similar to the one obtained without the
discrepancy. When including the discrepancy, the error decreases to $e = 0.0185$, which is not as
small as that obtained without the use of a surrogate (it is the price to pay to reduce
computational costs), but is still smaller than any error obtained without the inclusion of the
discrepancy. The use of a surrogate also impacts both the uncertainty of the parameters (the difference between the predictions using mean and MAP) and the predictions using the model and the discrepancy (broader band). The prediction using the surrogate yields an error of $e = 0.0116$, which is larger than the error obtained without the inclusion of the discrepancy, although we should recall that the use of a surrogate model here is merely instrumental.

\begin{figure}[tb]
     \centering
     \begin{subfigure}[b]{0.48\textwidth}
         \centering
         \includegraphics[width=\textwidth]{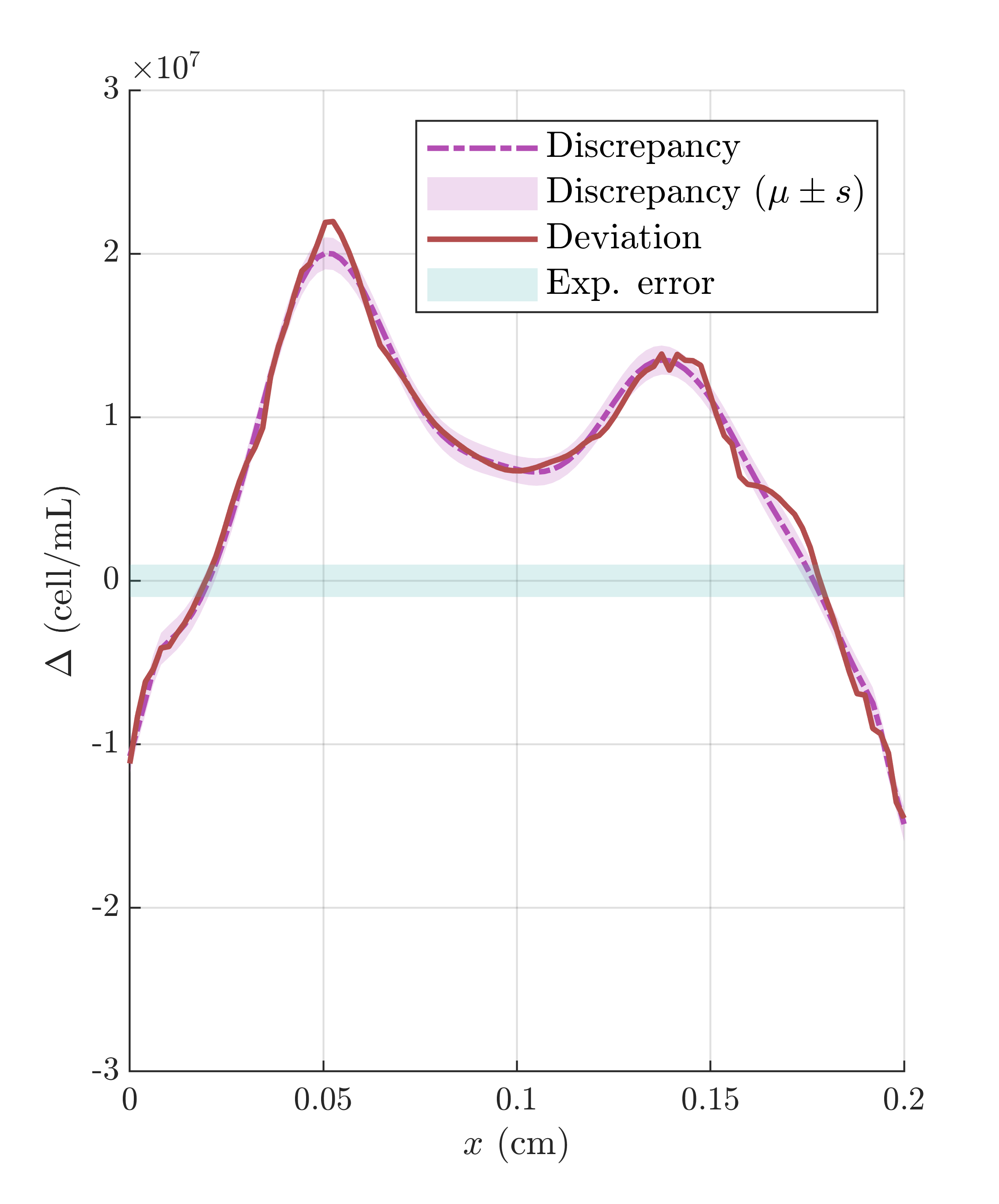}
         \caption{BCD approach.}
         \label{fig:inadequacy_B}
     \end{subfigure}
     \hfill
     \begin{subfigure}[b]{0.48\textwidth}
         \centering
         \includegraphics[width=\textwidth]{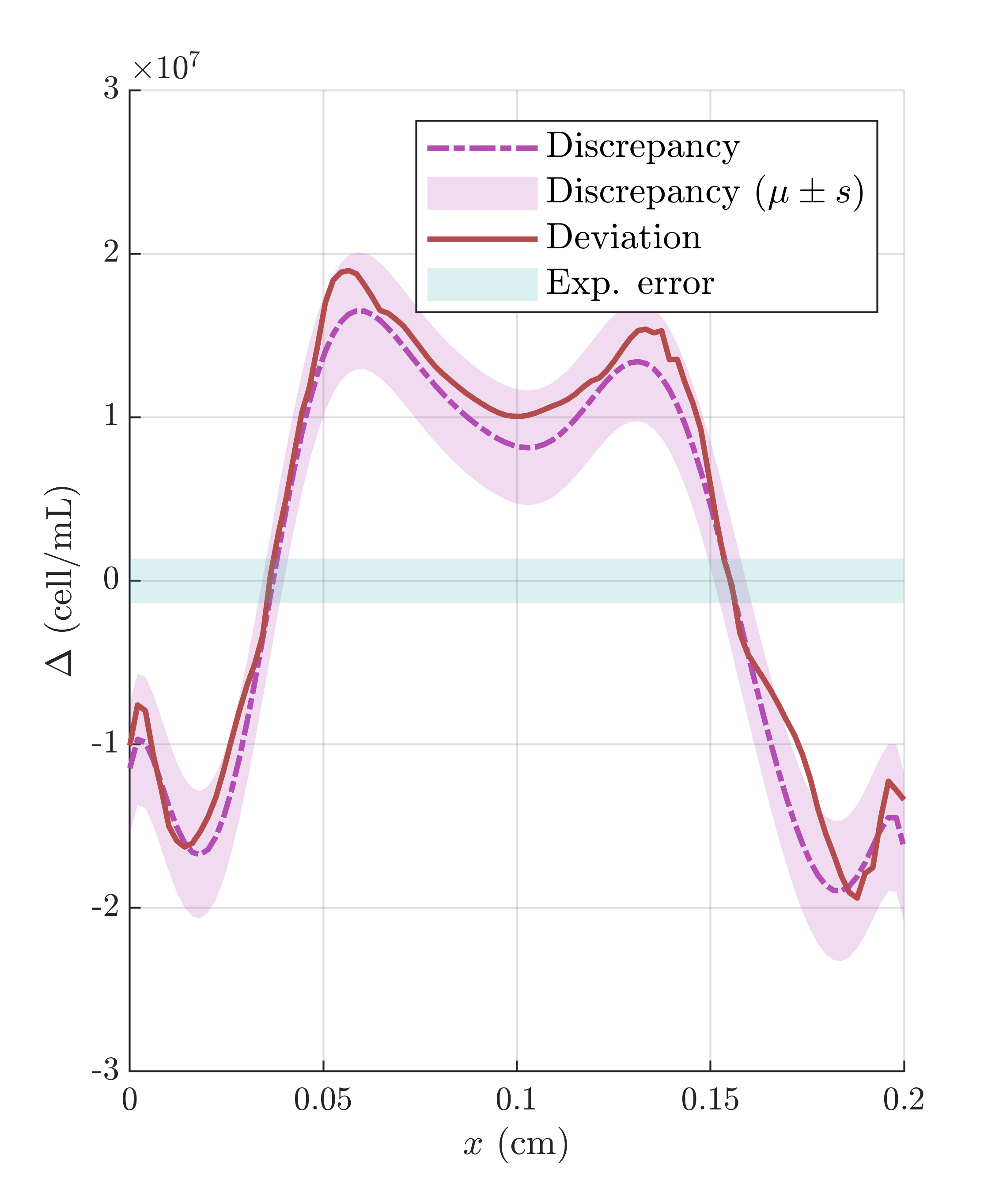}
         \caption{BCED approach.}
         \label{fig:inadequacy_D}
     \end{subfigure}
     \hfill
        \caption{\textbf{Evaluation of the model adequacy.} The discrepancy indicates the model inability to reproduce the experimental data; the deviation indicates the difference between the calibrated model and the data.}
        \label{fig:inadequacy}
\end{figure}

As discussed previously, the inclusion of a discrepancy error $\delta$ accounts for
the inability of the parametric model $\eta$ to reproduce the actual phenomenon, even with the best
possible choice of parameters. This allows for an initial assessment of the suitability of a model without the need for validation scenarios. In Fig. \ref{fig:inadequacy}, we illustrate this fact by showing the GP representing the discrepancy $\delta(x)$, the standard deviation of the experimental error $\epsilon(x)$, and the deviation $d(x) = z(x) - \eta(x;\hat{\boldsymbol{\theta}})$, where $\hat{\boldsymbol{\theta}}$ are the MAP estimation of the parameters. In both approaches, BCD and BCED, the discrepancy is close to the deviation, showing a good predictive capacity of the overall approach, whereas it is far from the experimental error, illustrating the limited ability of the parametric model to explain the data. In other words, we can recreate the process $\varphi$ with good accuracy, but there is room for improvement for an interpretable or mechanistic description. In fact, the use of $\delta$ instead of $d$ is more enlightening, since $d$ also incorporates the measurement error $\epsilon$. Decision-making is less obvious when using a surrogate, which is the price to pay for reducing computational cost.

In particular, for the case study analyzed, we observe negative bias close to the boundaries of the microfluidic device and positive bias at the center of the chamber, hinting at mechanisms that could be modified or enhanced in the parametric model.
\subsection{Using the probabilistic model for downstream tasks}
Once the parametric model is calibrated, the posterior distribution of the model parameters $\boldsymbol{\theta}$, denoted by $\mathcal{F}_{\boldsymbol{\theta}}$, may be used for many downstream tasks, such as uncertainty quantification of a certain Quantity of Interest (QoI) (velocity of a  propagating wave, GBM tumor burden growth, amount of cells extravasating the lateral channels, that is an indicator of the metastasis process, among others), goodness of fit evaluation of statistical tests (for instance to evaluate if a certain drug affects cell behavior through some model parameter). Therefore, an appropriate characterization of the statistics of the posterior distribution of the parameters is paramount. Here we illustrate this issue by analyzing the distributions $\mathcal{F}_{\boldsymbol{\theta}}$ for each of the approaches, as well as the distribution of the experimental error $\sigma$ and the eventual hyperparameters. We also discuss the statistics of the cell profile $\hat{u}(x), x \in [0;L]$, $L = 2 \, \mathrm{cm}$, as an example of QoI.

\begin{figure}[tb]
    \centering
    \includegraphics[width=\linewidth]{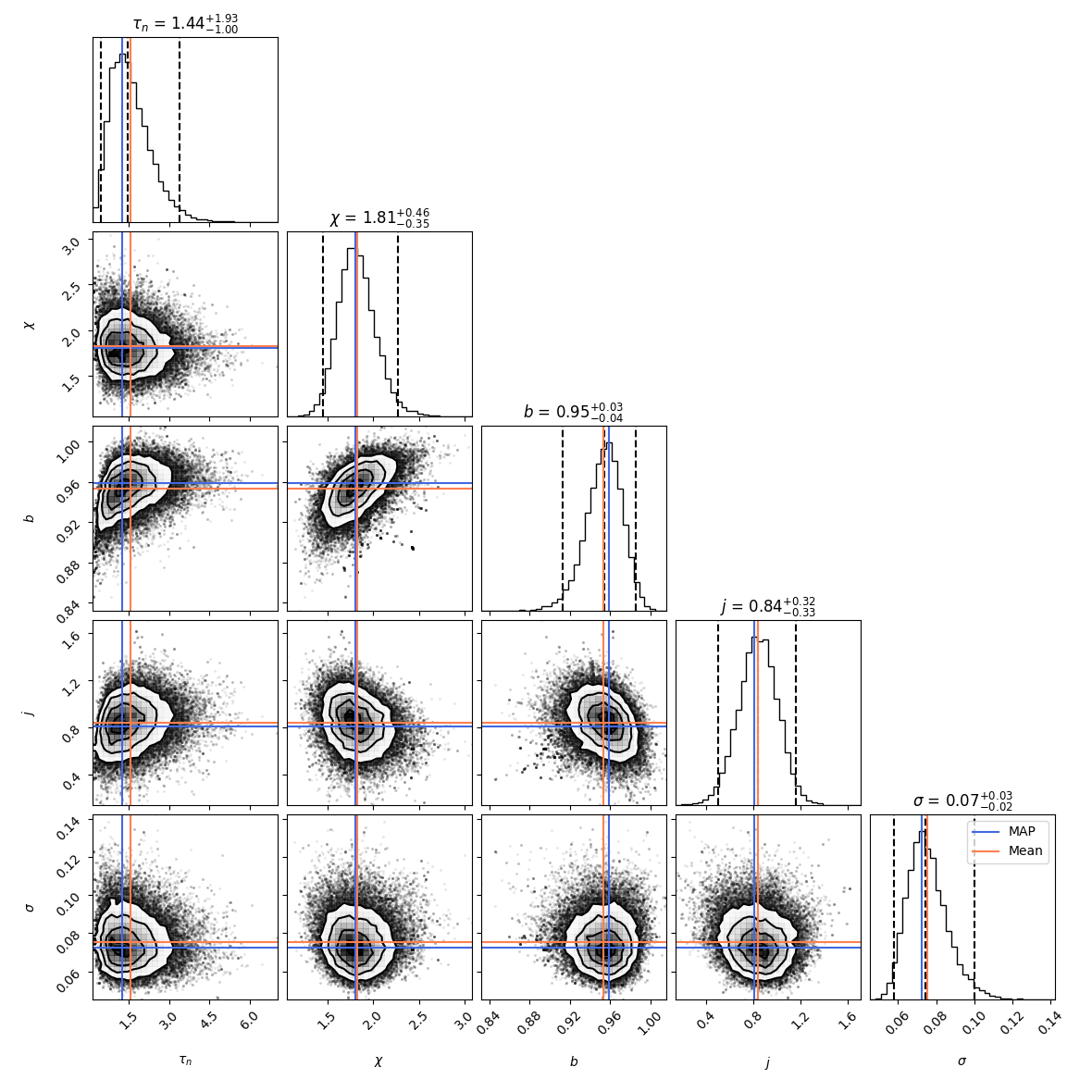}
     \caption{\textbf{Posterior distribution for the BI approach.} The posterior distributions are obtained for 8000 samples, 16 chains, and 0.2 burn-in.}
    \label{fig:corner_typeA}
\end{figure}

The results of the BI approach are summarized in the corner plot shown in Figure~\ref{fig:corner_typeA}. The first four parameters ($\tau_n, \chi, b, j$) correspond to the model parameters, while $\sigma$ quantifies the standard deviation of the experimental error. On the diagonal, the marginal posterior distributions for each parameter are displayed, with the median value and corresponding 95\% credible intervals indicated. The off-diagonal subplots illustrate the joint posterior distributions, capturing potential correlations between parameter pairs. Each of the model parameters is reasonably well constrained by the data, with unimodal marginal posteriors and moderate credible intervals. The standard deviation of the experimental error, $\sigma$, has a relatively narrow distribution, suggesting a robust estimation of the measurement noise and indicating that the model is capable of explaining the majority of experimental variability. Inspecting the joint distributions, weak to moderate correlations are apparent between some parameter pairs, e.g., between $b$ and $j$ or $\chi$ and $b$, while other pairs, such as $j$ and $\sigma$, show minimal dependency, which supports the identifiability of the parameter set under the current data and model configuration. MAP and mean estimates for all parameters are marked in blue and red, and found to be consistent with each other, further reflecting the unimodal nature of the posterior distributions, indicating convergence and well-defined parameter estimates. Taken together, these results show that this BI approach delivers precise parameter estimates along with a transparent uncertainty quantification, capturing both model and experimental errors, thus improving the reliability of model predictions.

\begin{figure}[tb]
    \centering
    \includegraphics[width=\linewidth]{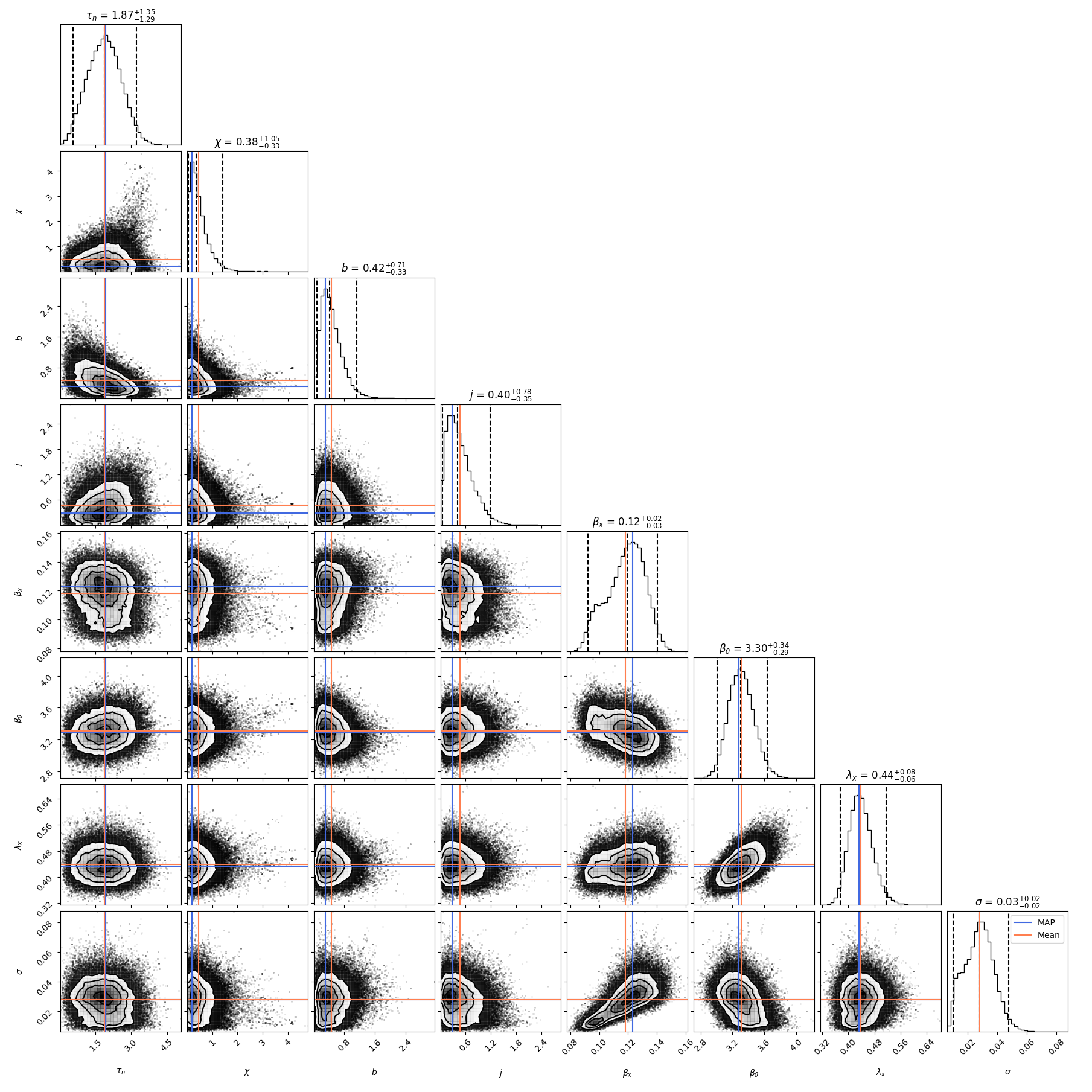}
    \caption{\textbf{Posterior distribution for the BCD approach.} The posterior distributions are obtained for 20000 samples, 16 chains, and 0.2 burn-in.}
    \label{fig:corner_typeB}
\end{figure}

The corner plot in Figure~\ref{fig:corner_typeB} presents the results of the BCE approach that include the first four parameters of the physical model ($\tau_n, \chi, b, j$), followed by three surrogate hyperparameters from the GP ($\beta_x, \beta_{\theta}, \lambda_x$), and finally the standard deviation of the experimental error, $\sigma$.
The physical model parameters are well constrained with unimodal marginal posteriors, indicating the data provide strong information about the system behavior.
The GP hyperparameters, which govern the surrogate model’s kernel structure and smoothness, exhibit distinct posterior distributions reflecting uncertainty in surrogate model tuning. Some weak to moderate correlations appear mainly among the hyperparameters and between hyperparameters and some model parameters, visible as elliptical shapes in the joint distributions. This reveals interdependencies in how surrogate model tuning interacts with physical parameters.
The experimental error standard deviation, $\sigma$, is tightly estimated, indicating a clear distinction between measurement noise and model-based uncertainty. Consistent MAP (blue) and mean (red) values across all parameters underscore unimodality and stable convergence of the BCE approach.

\begin{figure}[tb]
    \centering
    \includegraphics[width=\linewidth]{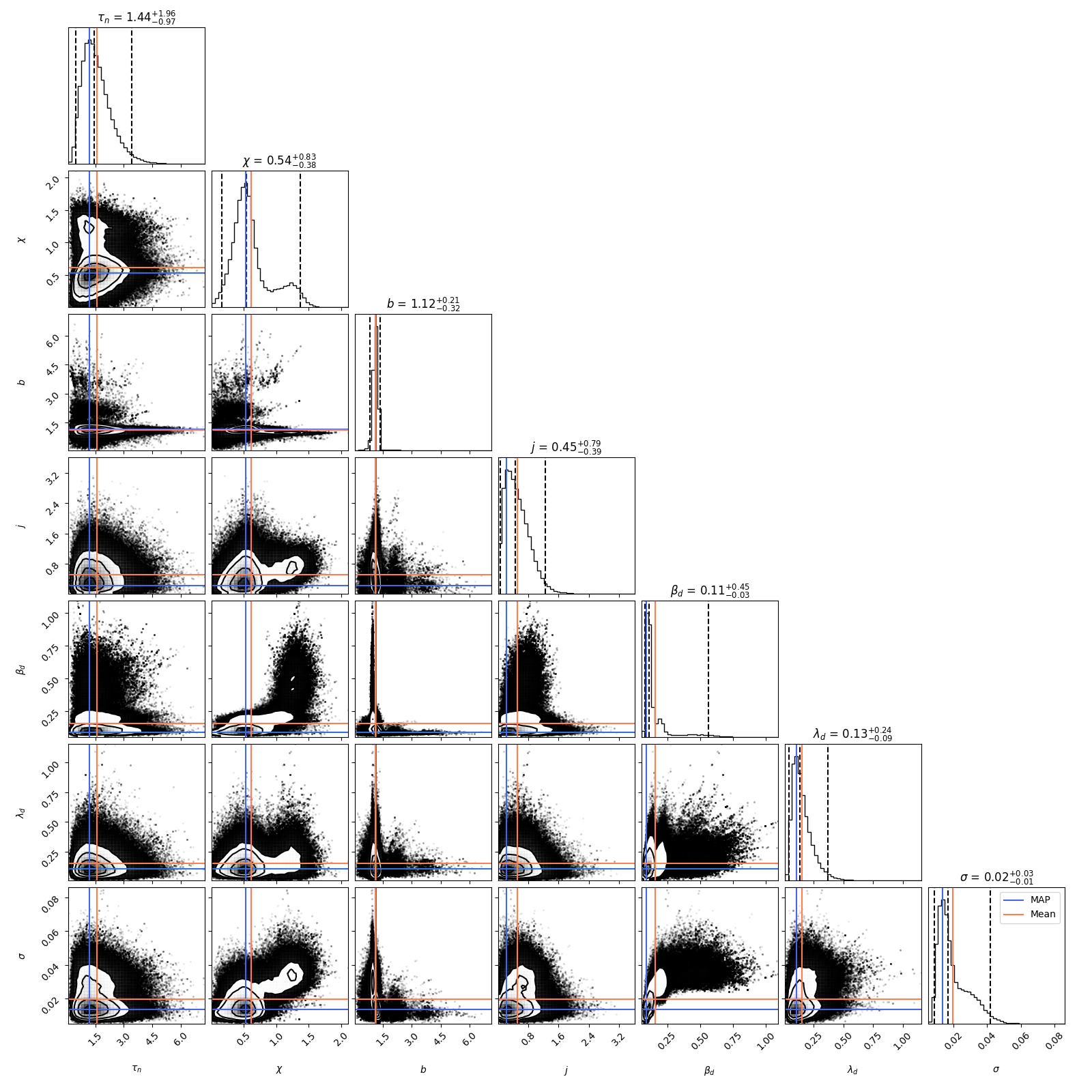}
    \caption{\textbf{Posterior distribution for the BCE approach.} The posterior distributions are obtained for 100000 samples, 16 chains, and 0.2 burn-in.}
    \label{fig:corner_typeC}
\end{figure}

Figure~\ref{fig:corner_typeC} displays the results of the BCD approach, now extended to estimate the model discrepancy. The first four parameters ($\tau_n, \chi, b, j$) correspond to the physical model, while the next two ($\beta_d, \lambda_d$) are the GP hyperparameters of the model discrepancy, and as before, $\sigma$ denotes the standard deviation of the experimental error.
The posterior distributions for $\tau_n, \chi, b,$ and $j$ remain unimodal; however, the credible intervals are generally wider compared to schemes without explicit discrepancy modeling. This reflects the additional uncertainty introduced by acknowledging a possible mismatch between model and reality.
The GP hyperparameters for model discrepancy, $\beta_d$ and $\lambda_d$, display broad and asymmetric posterior distributions, indicating significant uncertainty about the nature and scale of the discrepancy. Their inference enables the calibration framework to more flexibly account for systematic model errors. The experimental noise standard deviation, $\sigma$, stabilizes at a relatively low value, suggesting that the remaining variability attributed to measurement errors is well separated from both, parameter and model discrepancy uncertainties. Consistency between MAP (blue) and mean (red) estimates across most parameters continues to support stable convergence and essentially unimodal posteriors, despite the added flexibility of the discrepancy term.
Several off-diagonal panels reveal correlations, especially between physical model parameters and the discrepancy hyperparameters, highlighting the statistical interplay between parameter estimation and model inadequacy correction.
Including an explicit model discrepancy component in the Bayesian calibration thus increases the overall predictive uncertainty, more transparently combining the model uncertainty with the potential inadequacy of the model structure. This leads to more realistic uncertainty quantities and improves the trustworthiness of the predictive inferences.

\begin{figure}[tb]
    \centering
    \includegraphics[width=\linewidth]{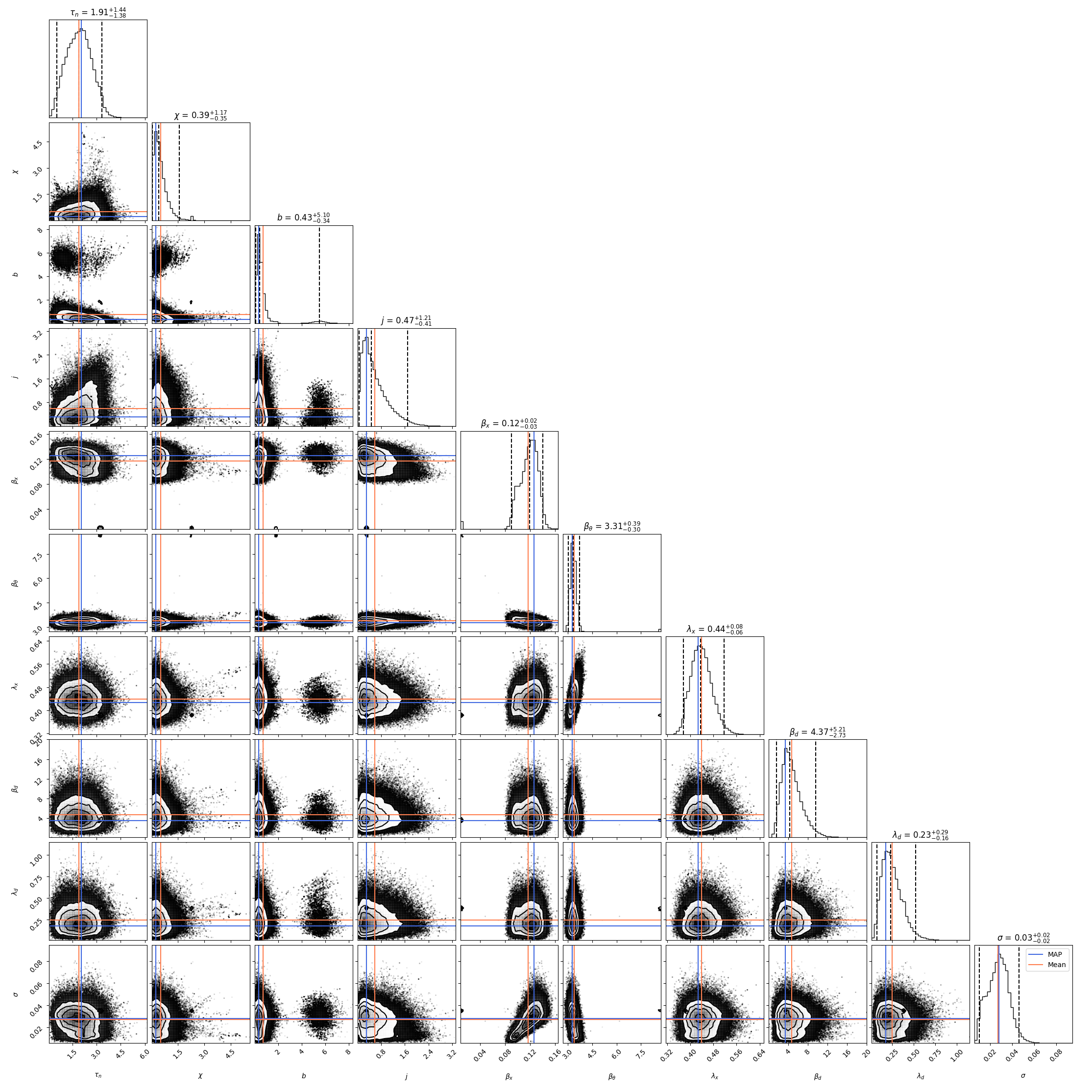}
    \caption{\textbf{Posterior distribution for the BCED approach.} The posterior distributions are obtained for \texttt{30000} samples, 32 chains, and 0.2 burn-in.}
    \label{fig:corner_typeD}
\end{figure}

Figure~\ref{fig:corner_typeD} illustrates the results of the BCED approach, incorporating both a
surrogate model and a model discrepancy term. The first four parameters ($\tau_n, \chi, b,
j$) come from the physical model; the next three ($\beta_x, \beta_{\theta}, \lambda_x$) are the GP
hyperparameters for the surrogate model, followed by two GP hyperparameters ($\beta_d, \lambda_d)$
for the model discrepancy; finally, the experimental error standard deviation, $\sigma$.
The marginal posterior distributions for all parameters are generally unimodal, indicating the parameters remain identifiable even with increased model flexibility.
The credible intervals for the physical and surrogate model parameters are somewhat broader compared to models without a discrepancy component, reflecting the honest propagation of combined uncertainties.
The GP hyperparameters for the surrogate ($\beta_x, \beta_{\theta}, \lambda_x$) and the discrepancy model ($\beta_d, \lambda_d)$ possess clearly informative posteriors, demonstrating that the data provide nontrivial constraints on not just the physical and surrogate models but also the nature and scale of the discrepancy.
Notable correlations are present between some surrogate and discrepancy hyperparameters, as evidenced by the elliptical shapes in their joint posterior panels. These dependencies illustrate the statistical interaction among model fidelity, surrogate expressiveness, and the allocated flexibility to correct for model inadequacy.
The standard deviation of the experimental error, $\sigma$, remains tightly constrained and distinct from other sources of uncertainty, facilitating a clear partitioning between measurement noise and modeling uncertainty.
The MAP (blue) and mean (red) estimates are close for all parameters, supporting the unimodal nature of the posteriors.
In summary, the inclusion of both a surrogate and a model discrepancy term leads to a richer and more transparent quantification of predictive uncertainty, pointing in particular to potential shortcomings in the physical model and the surrogate representation.
\begin{figure}[p]
     \centering
     \begin{subfigure}[b]{0.48\textwidth}
         \centering
         \includegraphics[width=\textwidth]{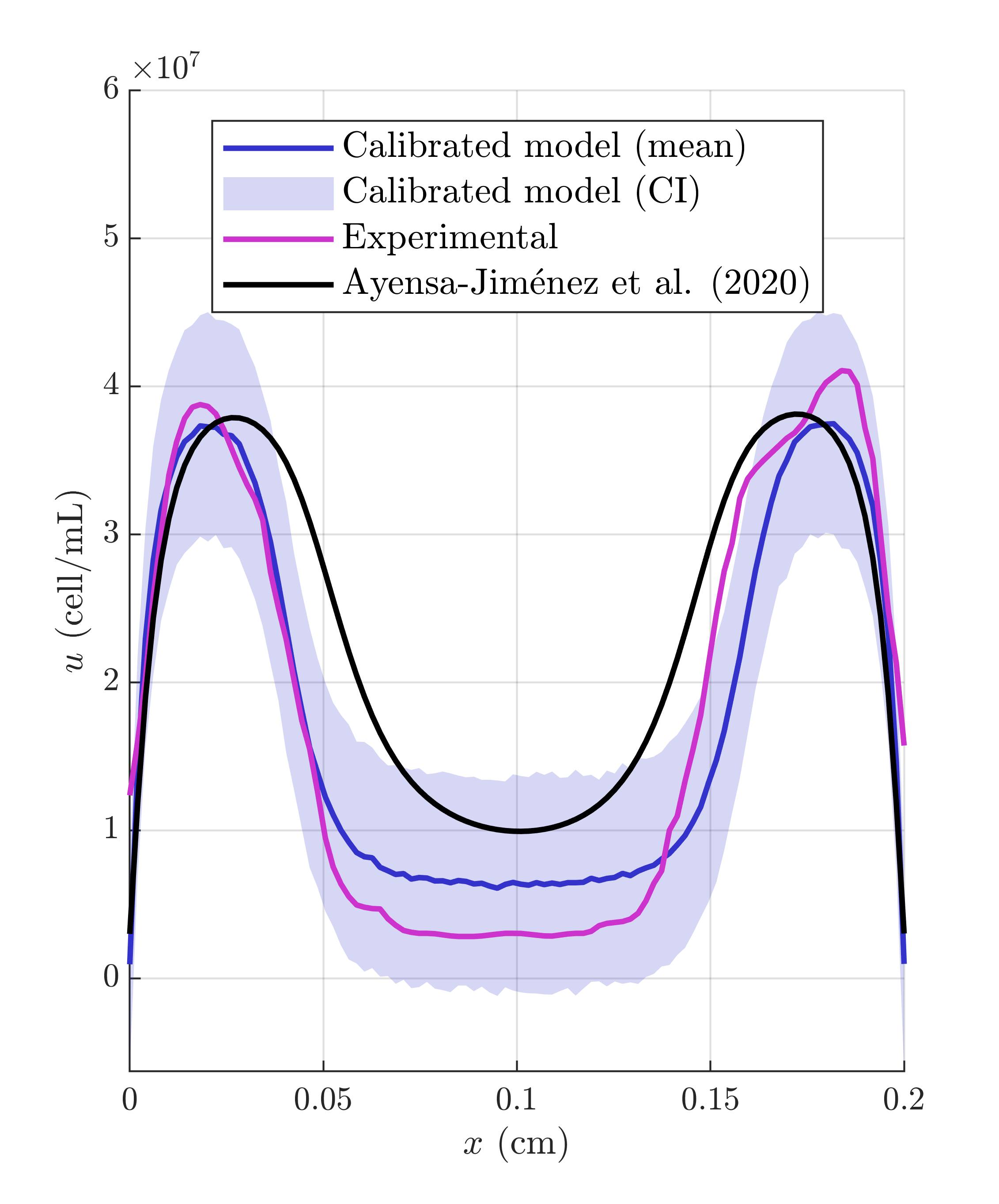}
         \caption{Parametric model.}
         \label{fig:sampling_A}
     \end{subfigure}
     \hfill
     \begin{subfigure}[b]{0.48\textwidth}
         \centering
         \includegraphics[width=\textwidth]{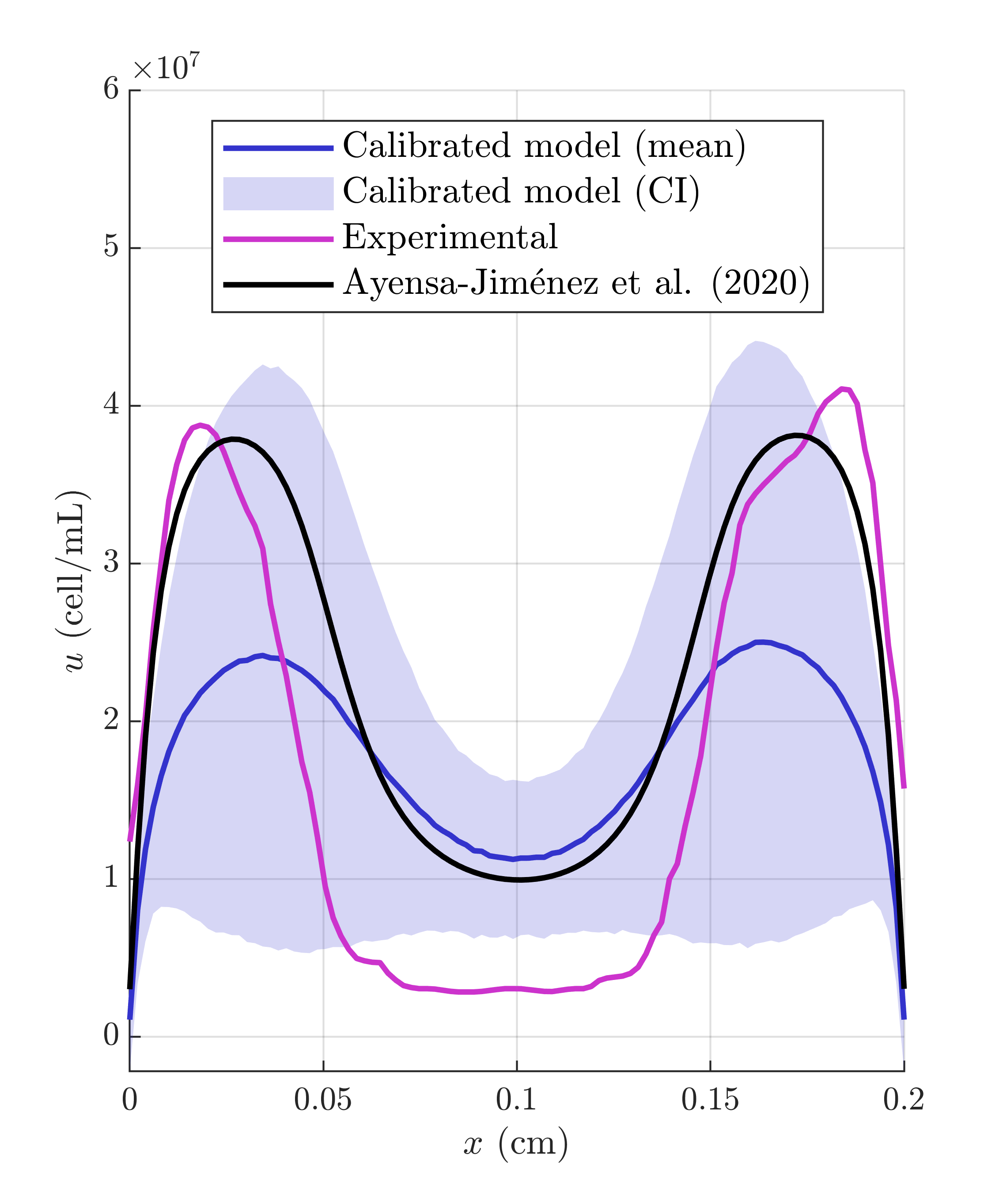}
         \caption{Parametric model and surrogate.}
         \label{fig:sampling_B}
     \end{subfigure}
     \hfill
     \begin{subfigure}[b]{0.48\textwidth}
         \centering
         \includegraphics[width=\textwidth]{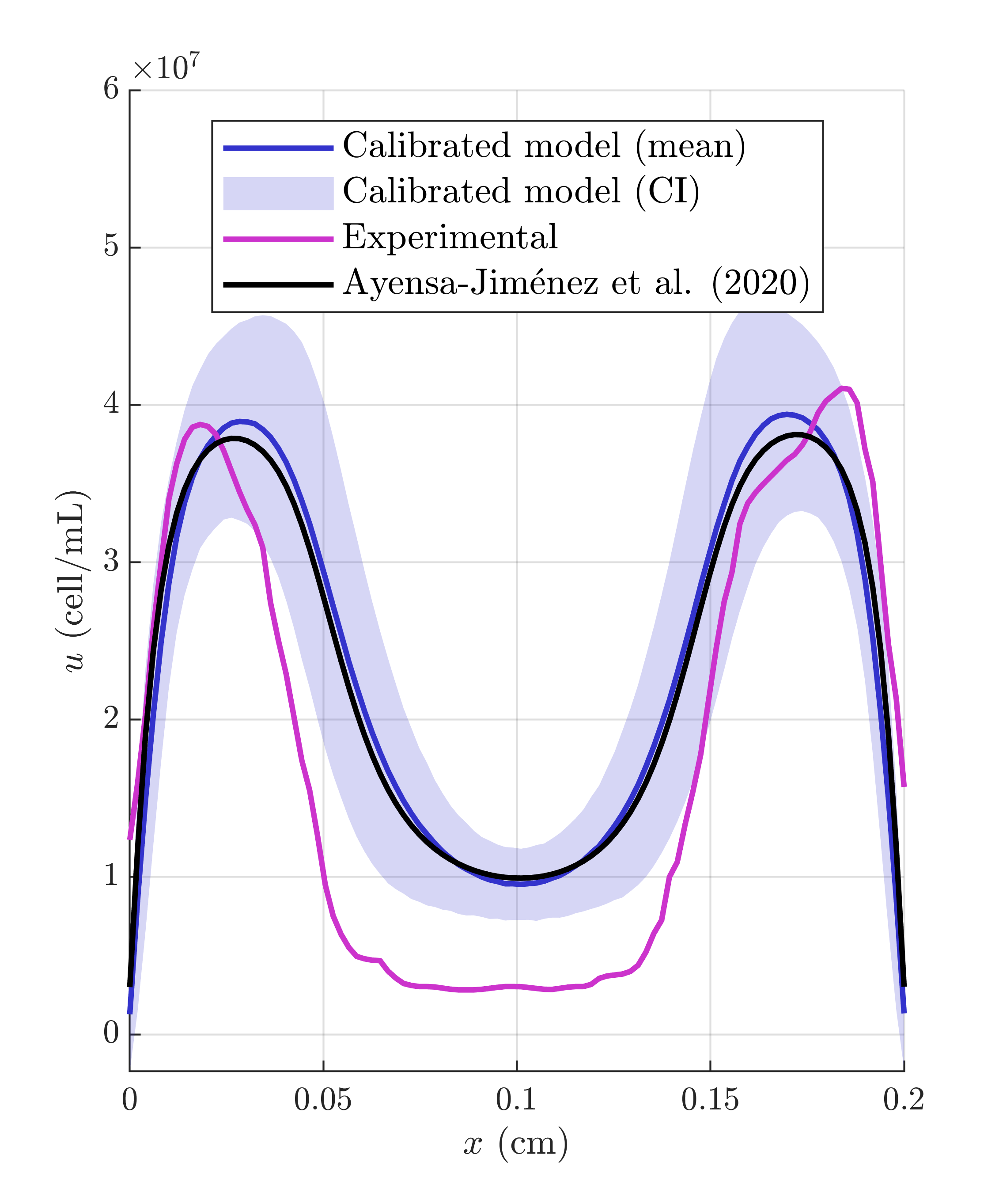}
         \caption{Parametric model and discrepancy.}
         \label{fig:sampling_C}
     \end{subfigure}
     \hfill
     \begin{subfigure}[b]{0.48\textwidth}
         \centering
         \includegraphics[width=\textwidth]{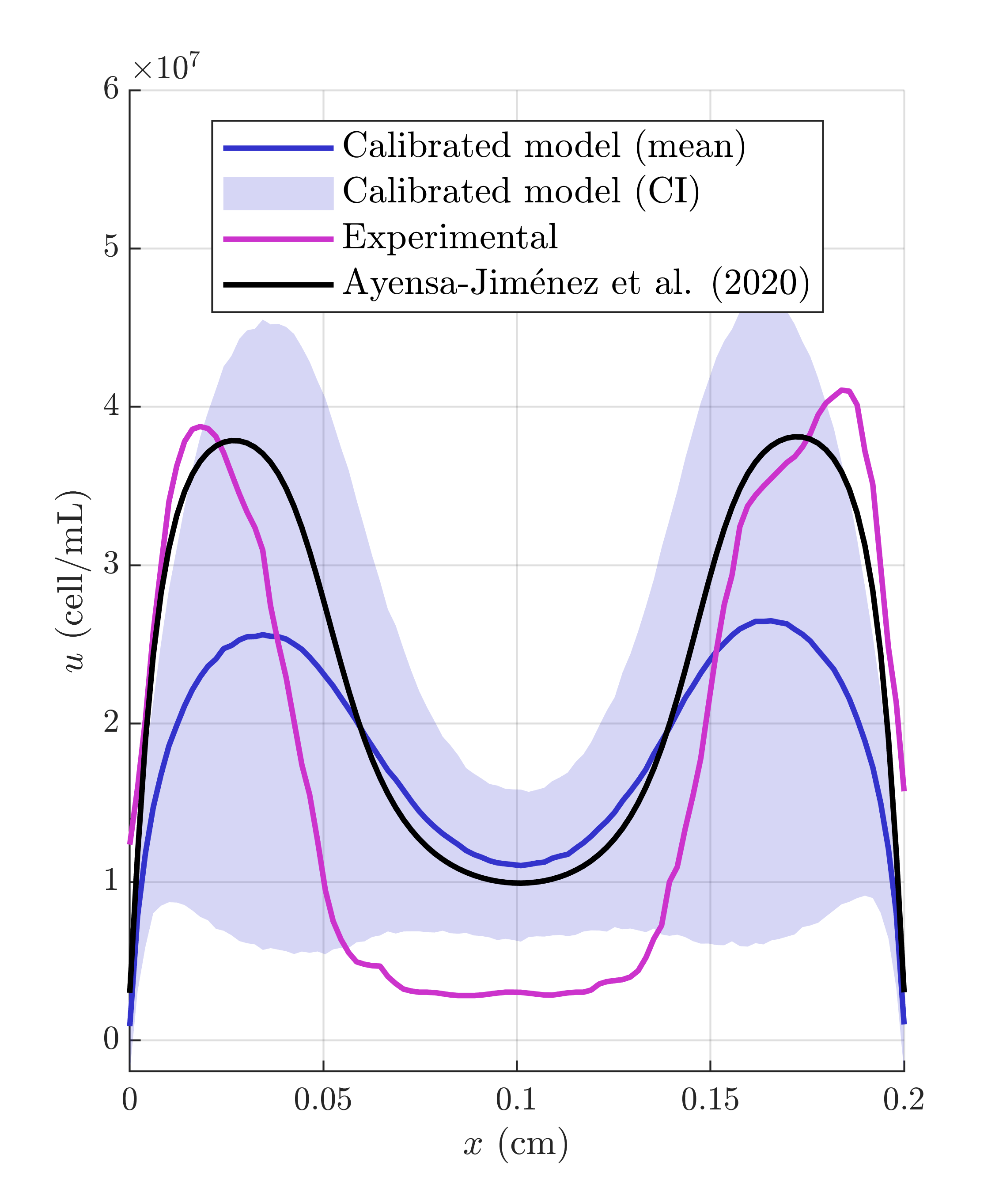}
         \caption{Par. model, surrogate and discrepancy.}
         \label{fig:sampling_D}
     \end{subfigure}
        \caption{\textbf{Robustness of the approach using the parametric model.} The robustness is evaluated by means of the prediction bias and the prediction uncertainties when considering the \emph{a posteriori} parameter distribution. }
        \label{fig:sampling}
\end{figure}

Including the discrepancy term has many practical consequences beyond a more appropriate characterization of the distribution of the parameters. In Fig. \ref{fig:sampling}, we show the cell profile and the experimental error obtained by Monte Carlo sub-sampling from the posterior distribution of the model parameters. More specifically, we generate $N$ cell profiles $u_i$, $i=1,\ldots,N$, where each profile is obtained by means of
\[u_i(x) = \hat{u}_i(x) + \epsilon_i,\]
where $\hat{u}_i(x) = u(x,\boldsymbol{\theta}_i)$ is the cell profile predicted using the parameter
$\boldsymbol{\theta}_i$ sampled from the posterior parameter distribution
$\mathcal{F}_{\boldsymbol{\theta}}$ and $\epsilon_i$ is sampled from the error distribution
$\mathcal{N}(0,\sigma^2)$. The figure shows the mean $\bar{u}(x)$ of a sample of size $N=1000$ as well as the region $\bar{u}(x) \pm 2s_u(x)$, where $s_u(x)$ is the standard deviation:
\[\bar{u}(x) = \frac{1}{N}\sum_{i=1}^N\hat{u}_i(x), \quad 
s_u(x) = \left(\frac{1}{N-1}\sum_{i=1}^N\left(u_i(x) - \bar{u}(x)\right)^2\right)^{1/2}.\]
\begin{figure}[htb]
     \centering
     \begin{subfigure}[b]{0.48\textwidth}
         \centering
         \includegraphics[width=\textwidth]{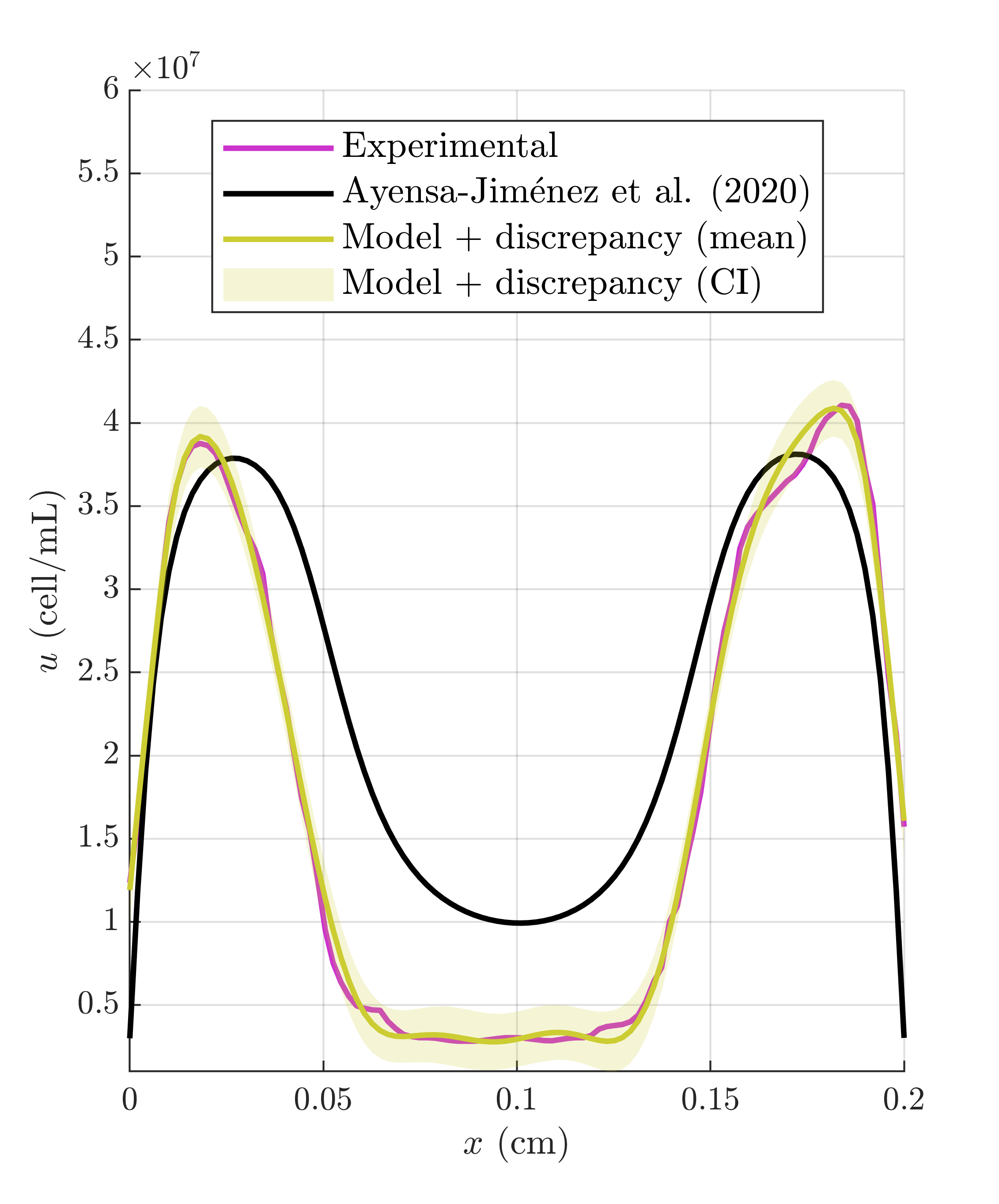}
         \caption{BCD approach.}
         \label{fig:discrepancy_C}
     \end{subfigure}
     \begin{subfigure}[b]{0.48\textwidth}
         \centering
         \includegraphics[width=\textwidth]{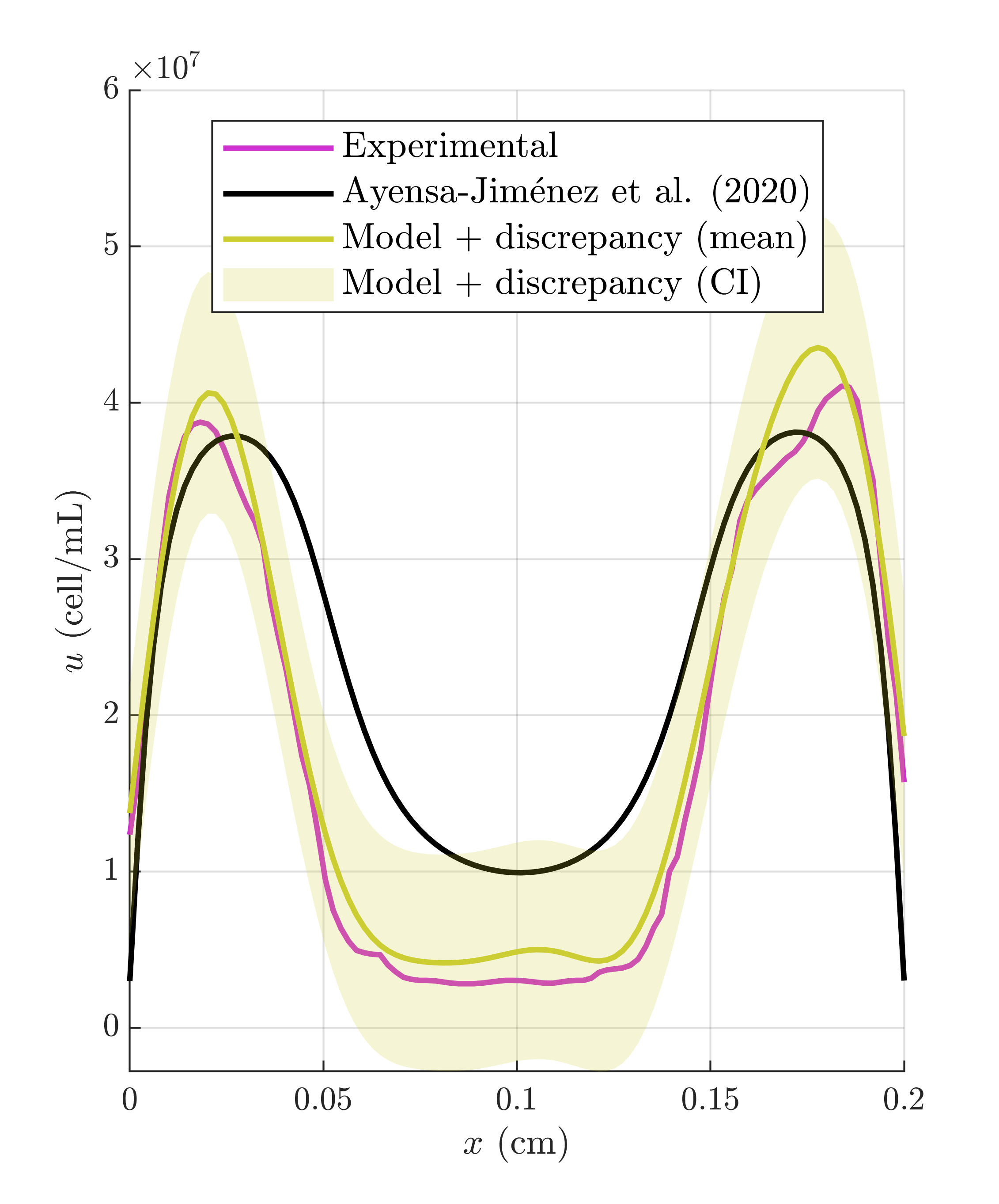}
         \caption{BCED approach.}
         \label{fig:discrepancy_D}
     \end{subfigure}
        \caption{\textbf{Prediction capacity when considering model discrepancy.} The robustness is evaluated by means of the prediction bias and the prediction uncertainties when considering the \emph{a posteriori} parameter distribution. }
        \label{fig:discrepancy}
\end{figure}

The benefit of including the discrepancy, in terms of uncertainty quantification, is illustrated in Fig.~\ref{fig:discrepancy}. In terms of predictive capacity, this approach can approximate the experimental data while keeping the uncertainty small. In fact, compared to the approach using the surrogate and the one using the inclusion of the discrepancy, we observe that the mean behavior is better described, and the confidence band defined as $\bar{u}(x) \pm 2s_u(x)$ is narrower and more regular.
\begin{figure}[p]
     \centering
     \begin{subfigure}[b]{0.48\textwidth}
         \centering
         \includegraphics[width=\textwidth]{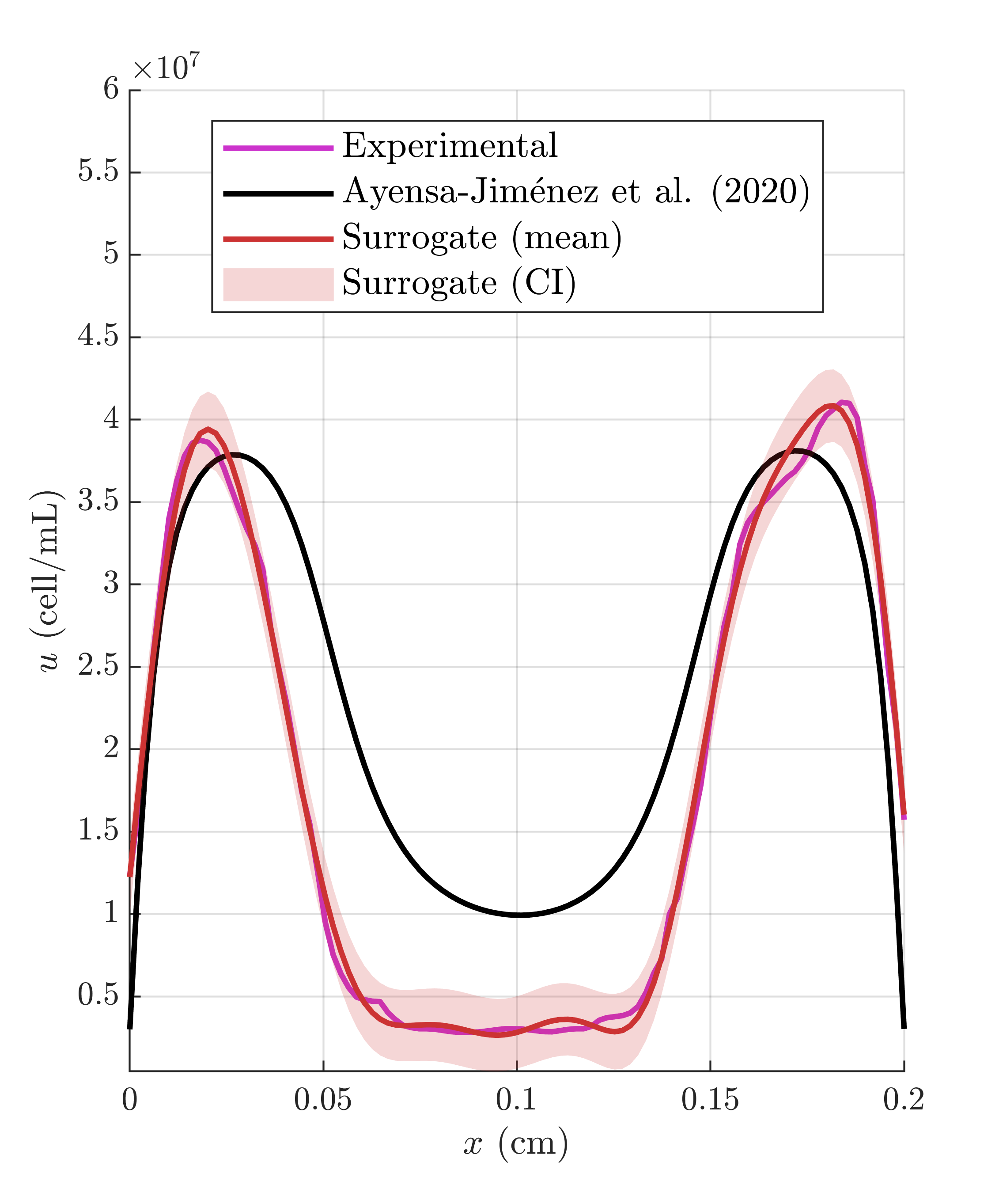}
         \caption{BCE approach.}
         \label{fig:surrogate_B}
     \end{subfigure}
     \begin{subfigure}[b]{0.48\textwidth}
         \centering
         \includegraphics[width=\textwidth]{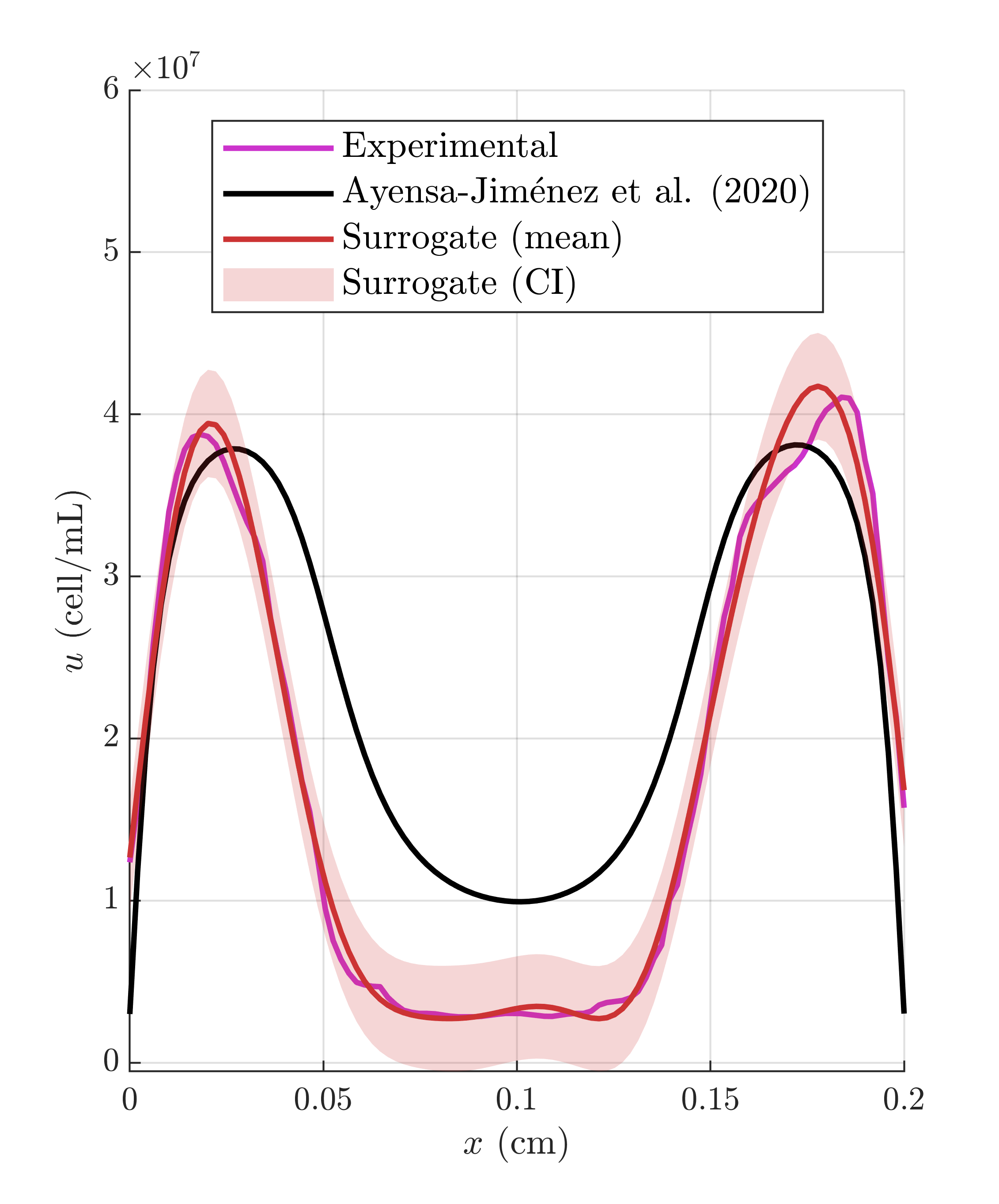}
         \caption{BCED approach.}
         \label{fig:surrogate_D}
     \end{subfigure}
        \caption{\textbf{Performance of the surrogate evaluation.} We show the ability of the GP for accurately predict the experimental data for the BCD and BCED approaches.}
        \label{fig:surrogate}
\end{figure}

As discussed before, in the presented approaches, the use of a surrogate model is merely instrumental to reduce computational cost. However, for completeness, in Fig. \ref{fig:surrogate} we show the ability of the GP to recreate the biological process $\varphi$. The results are very similar in both cases, but adding the discrepancy term slightly increases the bias and uncertainty because it lets the model adjust for differences between the surrogate and the true biological process. This correction adds flexibility but introduces uncertainty and shifts the predictions since the calibration must fit data and account for possible model errors, instead of fitting data alone. In other words, modeling discrepancy acknowledges imperfection in the surrogate, giving more reliable but less certain and slightly adjusted predictions.

The computational times of the studied algorithms were evaluated to provide a comparative perspective, acknowledging variability due to system specifications, network traffic, and other external factors. Specifically, BI required approximately 5 seconds per iteration; BCE achieved around 42 iterations per second, leading to a notable speedup using a surrogate model. BCD averaged around 5 seconds per iteration, while BCED reached approximately 22 iterations per second, also resulting in significant speed improvements through surrogate modeling. Compared to BI, the BCE approach achieves approximately a 208-fold speedup, while BCED attains about an 111-fold speedup relative to BCD. These improvements highlight the effectiveness of surrogate modeling in accelerating iteration performance. 

The discrepancy in our modeling framework provides valuable information about the performance and limitations of the physical model. Specifically, a large discrepancy highlights that the current physical model is not sufficiently accurate, indicating room for improvement in its formulation. Although incorporating discrepancy enriches our understanding of model inadequacies, it does not lead to improved estimates of the model parameters themselves, suggesting challenges such as parameter identifiability or model structural errors that cannot be resolved by calibration alone. Indeed, in Table \ref{table:parameters} we show the values of the parameters obtained using the BI and BCD approaches and we compare them to the ones reported in \cite{ayensa2020mathematical}. We focus on the parameters related with tumor progression ($\tau_n,\chi$ and $b$) and not with the experimental setup ($j$). Even if the BI approach shows the smallest error between the prediction given by the model and the experimental profile, the value of the chemotaxis parameter $\chi$ is overestimated. This is because the work in \cite{ayensa2020mathematical} employed a different experimental scenario (pseudopalisade configuration), for which GBM progression has more migratory characteristics, and for which the identifiability of $\chi$ was clearer. The inclusion of the discrepancy does not guarantee increased accuracy in the the estimated parameter, but it does indicate that the model calibration may not be appropriate, as here, without the need for a different validation scenario. The value of the remaining parameters remain close to their reference. It is interesting to analyze the results in light of previous studies that have focused on characterizing the  parameters related to go-or-grow \cite{ayensa2025overview}. The correlation that was reported between the parameters $\tau_n$, $\chi$ and $b$ because of the lack of identifiability for such a complex model (see for instance \cite{ayensa2020analysis,perez2023massive}) is also here revealed. However, both the spread and correlations of these parameters are reduced when including the discrepancy term, indicating that part of this uncertainty was due to model inadequacy. The variability of each parameter is also reduced when compared to the ranges found in the specialized literature \cite{ayensa2020mathematical}.
\begin{table}[htb]
\centering
\resizebox{0.75\linewidth}{!}{
\begin{tabular}{lccc}
\hline
       & Values by \cite{ayensa2020mathematical}                                                       & BI                                                                & BCD                                                                  \\
       \hline
$\tau_n$ ($\mathrm{s}$)& $7.5 \times 10^5$                                   & $6.5 \times 10^5 $                                   & $6.1 \times 10^5 $                                     \\
$\chi$ ($\mathrm{cm^2 \cdot mmHg^{-1} \cdot s^{-1}}$) & $7.5 \times 10^{-9}$ & $20 \times 10^{-9}$ & $5.8 \times 10^{-9}$   \\
$b$ ($\mathrm{mmHg}^{-1}$)    & $0.14$                                       & $0.14$                                      & $0.16$                                         \\
$J$ ($\mathrm{s \cdot cm^{-1}}$)    & $1 \times 10^{6}$                      & $2.1 \times 10^{6}$                   & $0.6 \times 10^{6}$ \\
\hline
\end{tabular}}
\caption{\textbf{Values of the parameters estimated with and without the discrepancy inclusion.} The values reported are the MAP using physical units.}
\label{table:parameters}
\end{table}

Using a surrogate model offers an alternative, computationally cheaper, that can be more predictive than the physical model. However, this benefit comes at a cost: the surrogate may compromise the quality of parameter estimates since it does not explicitly embed the underlying physics but relies instead on information derived from synthetic data. Consequently, while the surrogate is useful for saving computational cost, it lacks the interpretability and physical grounds necessary for robust parameter inference.

By employing a Bayesian calibration approach that integrates different methodologies, we obtain a more complete and nuanced picture of the modeling problem. This analysis enables quantification of various sources of uncertainty, including those from the physical model, the surrogate, and observational noise. As a result, we gain insights into the strengths and weaknesses not only of the original physical model, but also of the surrogate, informing subsequent model development and uncertainty management strategies.
In practice, these findings emphasize the need to carefully select modeling approaches based on the specific goals of a study. Surrogate models are well-suited for prediction tasks that demand computational efficiency but lack interpretability, whereas physical models remain essential when physical interpretability and accurate parameter estimation are prioritized.

\section{Conclusion}\label{sec:conc}
This paper has investigated the interplay of various Bayesian calibration
strategies, including the use of parametric and surrogate models, with and without 
discrepancy, in the context of models for glioblastoma
progression.

We first describe the modeling framework and detailed the calibration methodologies employed. Through comparative analysis of prediction accuracy, parameter identifiability, and uncertainty quantification, our results reveal the distinct strengths and complementarities of each approach, providing a detailed understanding of how these modeling choices impact both practical outcomes and theoretical insights.

Our findings demonstrate that leveraging Bayesian inference for this glioblastoma progression case study can provide complementary perspectives: parametric models retain interpretability and physical grounding, while surrogates deliver computational efficiency and strong predictive performance. The explicit consideration of model discrepancy enriches the process by highlighting when model structure is insufficient and where improvements are needed. A key outcome of our investigation is the realization that these calibration strategies serve distinct but interconnected goals. While standard GP surrogate models excel in accurate and efficient prediction, particularly when computational resources are limited, they are less reliable for parameter inference due to the lack of embedded physical constraints. Conversely, calibration of the physical model, especially when accounting for discrepancy, provides greater insight into parameter identifiability and model limitations, revealing structural errors that cannot be addressed by calibration alone.

By systematically quantifying uncertainties from all sources -- parametric model, surrogate, and observational noise -- our Bayesian framework allows for a nuanced selection of modeling strategies based on whether prediction or physical understanding is prioritized. The results here illustrate that integrating multiple calibration methodologies not only improves predictive skill but also builds a deeper, more transparent understanding of model trustworthiness and limitations.

Overall, this work highlights the importance of an integrated calibration approach tailored to the specific objectives of the considered study. For future research on these types of glioblastoma progression models, we recommend continued development of hybrid and physics-informed surrogate models, as well as refinement of physical models to enhance their descriptive capacity and reduce structural errors. These efforts are also expected to foster calibration and prediction strategies that are less tightly bound to specific goals, promoting broader applicability across modeling scenarios.
\section*{Author contributions}
All authors have made substantial intellectual contributions to the study conception, execution, and design of the work. All authors have read and approved the final manuscript. In addition, the following contributions occurred: 
Conceptualization: Christina Schenk (C.S.), Ignacio Romero (I.R.); Methodology: C.S., Jacobo Ayensa-Jim\'enez (J.A-J.), I.R.; Formal analysis and investigation: C.S., J.A-J.; Writing: C.S., J.A-J., I.R.; original draft preparation: C.S., J.A-J., I.R.; Writing - review and editing: C.S., J.A-J., I.R.; Data Curation: C.S., J.A-J.; Visualization: C.S., J.A-J.; Software: C.S., J.A-J., I.R.; Funding acquisition: J.A-J..
\section*{Data and code availability statement}
The data and the Matlab code that support the findings of this study are openly available in the GBM-SIMUL GitHub repository \cite{ayensa2025} at \url{https://github.com/Ayensa-Jimenez/GBM-SIMUL}. The Bayesian calibration code is provided as an attached .zip file for reviewer access and will be made publicly available through an open repository upon publication of this work.
\section*{Acknowledgments}
The authors acknowledge financial support from the Spanish Ministry of Science and Innovation through a Juan de la Cierva Postdoctoral Fellowship for J.A-J. (Grant No. JDC2023-052319-I). The authors acknowledge the use of the large language model Perplexity AI to assist in reformulating the wording of several sentences in the manuscript for clarity. 

\clearpage
\bibliographystyle{siamplain}
\bibliography{refs}

\appendix
\section{Computational implementation}
GBM model was implemented in Matlab R2023b using the \texttt{pdepe} solver for systems of 1-D parabolic and elliptic PDEs. The \texttt{pdepe} routine uses a space-time integrator based on a piecewise nonlinear Galerkin approach that is second order accurate in space \cite{skeel1990method}, computing the solution over a specified spatial mesh of size $N_x = 100$ and time span with $N_t = 100$ (although the numerical solver is adaptive in time). Linear interpolation of the experimental data corresponding to the initial condition was performed using \texttt{interp1} to map the initial cell profile to the values at given spatial coordinates. The GBM progression simulation code is free and publicly available \cite{ayensa2025} (\url{https://github.com/Ayensa-Jimenez/GBM-SIMUL}).

The Bayesian calibration procedure was performed in Python using ACBICI 2.0.0, an in-house software package \cite{acbici}, leveraging core scientific libraries including \texttt{scipy}, \texttt{numpy}, and \texttt{math}. Experimental design for synthetic datasets was conducted using the \texttt{pyDOE} package, employing Latin Hypercube sampling to ensure efficient coverage of the parameter space. For Bayesian inference, the \texttt{emcee} package \cite{Foreman_Mackey_2013,Foreman-Mackey2019} was utilized to perform Markov Chain Monte Carlo (MCMC) sampling, enabling efficient exploration of the posterior distribution. Graphical illustrations and visualizations were generated with the \texttt{matplotlib}, \texttt{corner} and \texttt{seaborn} libraries.
ACBICI 2.0.0 and the code used for the case studies are provided in an attached .zip file for reviewer access and will be made publicly available in an open repository upon publication of this work.

\section{Data and model configurations}
\subsection{Data preparation}
A selection of experimental data points was performed to ensure comprehensive coverage of the input domain while maintaining computational efficiency. Initially, the data points corresponding to the minimum and maximum values of the key input variables were identified and included to anchor the dataset boundaries. Subsequently, 28 additional points were selected using Latin Hypercube Sampling (LHS) to uniformly cover the input space in a stratified manner.

For synthetic data generation, a set of 200 data points was selected at random from a larger pool of 500 simulations conducted in MATLAB. These synthetic simulations were performed by sampling the model parameters uniformly in the nominal interval $[0.1,6]$. This selection reflects plausible variations of the parameters within the identified operational range and provides a sufficiently large and diverse surrogate dataset compromising sufficient knowledge and computational efficiency for subsequent model calibration and validation.
\subsection{Priors}
The default priors from ACBICI 2.0.0 were chosen. Specifically, the hyperparameters were assigned $\mathrm{Gamma}(\alpha,\beta)$ priors scaled to the characteristic magnitudes of the data.
When no prior information was available, the shape parameter was fixed at $\alpha=5$, while the rate parameter $\beta$
was determined from empirical statistics of the inputs or outputs.
For the length-scale parameters $\beta_x$ and $\beta_t$, the rate was set to $\beta = 5/d$,
with $d$ denoting the average pairwise distance in the experimental inputs or synthetic time coordinates, respectively.
For the variance-related parameters $\lambda_x$ and $\lambda_d$, the rate was chosen as $\beta = 5/\sigma$,
where $\sigma$ is the empirical standard deviation of the experimental responses.

For the observational error standard deviation $\sigma_\varepsilon$, the prior was constructed from the mean absolute
response amplitude, $\bar{y} = \mathrm{mean}(|y^{\mathrm{exp}}|)$.
The prior mean was set to $0.1 \, \bar{y}$, with a standard deviation equal to $0.1$ of this prior mean.
These values defined the parameters of a Gamma distribution, ensuring that the prior expectation and
variance reflected a small but non-negligible fraction of the observed signal magnitude.
This specification yields weakly informative priors that are scaled to the natural variability of
the observed data, preventing unrealistic hyperparameter values while retaining calibration flexibility.
\end{document}